\documentclass[12pt]{amsart}

\usepackage{amsmath,amsfonts,amssymb}

\theoremstyle{plain}
\newtheorem{lemma}{Lemma}[section]
\newtheorem{theorem}[lemma]{Theorem}

\theoremstyle{definition}

\numberwithin{equation}{section}
\DeclareMathOperator{\Var}{Var}
\DeclareMathOperator{\var}{var}
\DeclareMathOperator{\diam}{diam}

\begin{document}

\newcommand{\mB}{\mathcal{B}}
\newcommand{\omB}{\overline{\mathcal{B}}}
\newcommand{\oL}{\overline{\Lambda}}
\newcommand{\graph}{(\oL,\omB(\Lambda))}
\newcommand{\ZZ}{\mathbb{Z}^{2}}
\newcommand{\lra}{\leftrightarrow}
\newcommand{\lrad}{\overset{*}{\leftrightarrow}}
\newcommand{\lrap}{\overset{+}{\leftrightarrow}}
\newcommand{\lram}{\overset{-}{\leftrightarrow}}

\title[Ising Spectral Gap]
{The Spectral Gap \\
of the 2-D stochastic Ising model \\ 
With Nearly Single-Spin Boundary Conditions}
\author{Kenneth S. Alexander}
\address{Department of Mathematics DRB 155\\
University of Southern California\\
Los Angeles, CA  90089-1113 USA}
\email{alexandr@math.usc.edu}
\thanks{Research supported by NSF grant DMS-9802368.}

\date{\today}

\begin{abstract}
We establish upper bounds for the spectral gap of the
stochastic Ising model at low temperature in an $N \times N$ box, 
with boundary conditions which are ``plus'' except for small regions
at the corners which are either free or ``minus.''  The spectral gap decreases
exponentially in the size of the corner regions, when these regions are of size at
least of order $\log N$.  This means that removing as few as $O(\log N)$ plus
spins from the corners produces a spectral gap far smaller than the
order $N^{-2}$ gap believed to hold under the all-plus boundary condition.
Our results are valid at all subcritical temperatures.
\end{abstract}
\maketitle

\section{Introduction and Main Theorem} \label{S:intro}
Let $\Lambda \subset \mathbb{Z}^{2}$ and let $\eta \in \{-1,0,1\}^{\partial
\Lambda}$.  Here
$\partial \Lambda = \{x \in \mathbb{Z}^{2} \backslash \Lambda: x \text{ is adjacent
to some site }$ $\text{in } \Lambda\}$.  The Hamiltonian for the Ising model on
$\Lambda$ with boundary condition $\eta$ is
\[
  H_{\Lambda,\eta}(\sigma) = - \sum_{\langle xy \rangle: x, y \in \Lambda}
  \sigma_{x}\sigma_{y} - \sum_{\langle xy \rangle: x \in \Lambda, y \in \partial
  \Lambda} \sigma_{x}\eta_{y}, \quad \sigma \in \{-1,1\}^{\Lambda},
\]
where the first sum is over unordered pairs of adjacent sites.  Let
$\mu = \mu_{\Lambda,\eta}^{\beta}$ denote the equilibrium measure 
when the inverse temperature is $\beta$:
\[
  \mu_{\Lambda,\eta}^{\beta}(\sigma) = (Z_{\Lambda,\eta}^{\beta})^{-1} e^{-\beta
    H_{\Lambda,\eta}(\sigma)},
\]
where $Z_{\Lambda,\eta}^{\beta}$ is the partition function.
Let $\Sigma_{\Lambda} =
\{-1,1\}^{\Lambda}$ be the configuration space, and let $\sigma_{\Lambda}$
denote a generic configuration.   (When convenient we write
``+'' and ``-'' in place of 1 and -1 for the two spins.)  

We consider the time evolution of the dynamic version of the model under Glauber
dynamics.  Let $\sigma^{x}$ denote the configuration given by $\sigma^{x}_{y} =
-\sigma_{y}$ for $y = x$, $\sigma^{x}_{y} = \sigma_{y}$ for $y \neq x$.
The flip rate at a site $x$ when the configuration is $\sigma$ is denoted
$c(x,\sigma)$ (notationally supressing its dependence on $\eta$.)
We assume that the flip rates are uniformly bounded:
\[
  0 < c_{0}^{\prime} \leq c(x,\sigma) \leq c_{0} \quad \text{for all } x 
    \text{ and } \sigma.
\]
We also make the usual assumptions that the flip rates are attractive and
translation invariant, satisfy detailed balance
and have finite range; see e.g. \cite{Li85} for full
descriptions of these properties.  The generator
$A = A_{\Lambda,\eta}^{\beta}$ of the corresponding Markov process is given by
\[
  (Af)(\sigma) = \sum_{x \in \Lambda} c(x,\sigma) \bigl( f(\sigma^{x}) 
    - f(\sigma) \bigr),
\]
and the Dirichlet form $\mathcal{D} = \mathcal{D}_{\Lambda,\eta}^{\beta}$ by
\[
  \mathcal{D}(f,g) = \langle f,Ag \rangle_{\mu},
\]
so that
\[
  \mathcal{D}(f,f) = - \frac{1}{2}\sum_{x \in \Lambda} \sum_{\sigma \in
    \Sigma_{\Lambda}} \mu(\sigma) c(x,\sigma) \bigl( f(\sigma^{x}) - 
    f(\sigma) \bigr)^{2},
\]
where $\mu = \mu_{\Lambda,\eta}^{\beta}$.
The spectral gap
$\Delta = \Delta(\Lambda,\eta,\beta)$ for such dynamics, which is the smallest
positive eigenvalue of $-A^{\beta}_{\Lambda,\eta}$, has the following
representation:
\begin{equation} \label{E:variational}
  \Delta(\Lambda,\eta,\beta) = \inf_{f \in L^{2}(\mu)} 
    - \frac{\mathcal{D}(f,f)}{\var_{\mu}(f)},
\end{equation}
where $\var_{\mu}(f)$ denotes the variance of $f$.  The gap $\Delta$ describes the
rate of exponential convergence in $L^{2}(\mu^{\beta}_{\Lambda,\eta})$ to
equilibrium, in the sense that for $S(\cdot)$ the semigroup generated by $A$ and
$\| \cdot \|_{\mu}$ the $L^{2}(\mu)$ norm, $\Delta$ is the largest constant such
that
\[
  \| S(t)f - \int f \, d\mu \|_{\mu} \leq \| f - \int f \, d\mu \|_{\mu}
    e^{-\Delta t} \quad 
    \text{for all } f \in L^{2}(\mu) \text{ and } t \geq 0.
\]

We say that two
configurations $\sigma,\sigma^{\prime} \in \Sigma_{\Lambda}$ are
\emph{adjacent}  if for some $x \in \Lambda$ we have $\sigma_{x} \neq
\sigma^{\prime}_{x}$ but
$\sigma_{y} = \sigma^{\prime}_{y}$ for all $y \neq x$.  For $S \subset
\Sigma_{\Lambda}$ define
\[
  \partial_{in} S = \{\sigma \in S: \sigma \text{ is adjacent to some site of }
  S^{c} \}. 
\]
Considering only indicator functions we obtain
\begin{equation} \label{E:indicators}
  \Delta(\Lambda,\eta,\beta) 
    \leq c_{0}|\Lambda| \inf_{D \subset \Sigma_{\Lambda}}
    \frac{\mu^{\beta}_{\Lambda,\eta}(\partial_{in}
    D)}{\mu^{\beta}_{\Lambda,\eta}(D)
    (1 - \mu^{\beta}_{\Lambda,\eta}(D))}
\end{equation}

Let $\tilde{\Lambda}_{N} = [-N,N]^{2}$ and $\Lambda_{N} = \tilde{\Lambda}_{N} \cap
\mathbb{Z}^{2}$.  
Let $l_{h}$ and $l_{v}$ denote the horizontal and vertical axes,
respectively, and consider the boundary condition $\eta^{k,\epsilon}$ given 
for $k \geq 0$ and $\epsilon \in \{0,-1\}$ by
\begin{equation} \label{E:bdrycond}
  \eta^{k,\epsilon}_{x} = 
  \begin{cases}
    1, &\text{if } d(x,l_{h}) \leq k \quad \text{or} \quad
      d(x,l_{v}) \leq k; \\
    \epsilon, &\text{otherwise}.
  \end{cases}
\end{equation}
Here $d(\cdot,\cdot)$ denotes Euclidean distance.
As a special case of results in \cite{HY} we have that for $\beta$ very large, 
for some $C, \lambda$ depending only on $\delta, \beta$
and $\epsilon$,
\begin{equation} \label{E:smalldelta}
  \Delta(\Lambda_{N},\eta^{\delta N,\epsilon},\beta) \leq C e^{-\lambda N}
    \quad      
    \text{for all } \delta \leq \frac{1}{2} \text{ and } N \geq 1.
\end{equation}
Here we will generalize this as follows.  Let $\beta_{c}$ denote the critical
inverse temperature of the Ising model on $\mathbb{Z}^{2}$.

\begin{theorem} \label{T:main}
  Let $\beta > \beta_{c}$.

  (i) For some $C, K, \lambda$ depending only on $\beta$, for all $N \geq 1$ and
  $k \geq 1$ satisfying $N - k \geq K \log N$,
  \begin{equation} \label{E:smalldelta2}
    \Delta(\Lambda_{N},\eta^{k,0},\beta) \leq C e^{-\lambda (N - k)}.
  \end{equation}

  (ii) For some $C, K \lambda$ depending only on $\beta$, for all $N \geq 1$ and 
  $k \geq 1$ satisfying $\min(k,N-k) \geq K \log N$,
  \begin{equation} \label{E:smalldelta3}
    \Delta(\Lambda_{N},\eta^{k,-1},\beta) \leq C e^{-\lambda \min(k,N-k)}.
  \end{equation}
\end{theorem}

It is easily verified that the spectral gap changes by at most a
constant when a single boundary spin is changed; in particular this can be applied
in comparing $k = 1$ to a completely free boundary.  Thus Theorem \ref{T:main}(i)
essentially includes exponential decay of the gap under the free boundary condition
at low temperatures, a result obtained by Thomas \cite{Th}.

Theorem \ref{T:main} generalizes (\ref{E:smalldelta}) to all $\delta
< 1$ and $\beta > \beta_{c}$,
and shows that the $L^{2}$ rate of convergence to equilibrium in the stochastic 
Ising model can be quite slow even when the boundary condition is 
overwhelmingly of a single spin, and otherwise free.  For the full ``plus''
boundary condition ($k = N$) at subcritical temperatures, it is believed \cite{FH}
that the spectral gap is of order
$N^{-2}$; Martinelli \cite{Ma} proved that for very low temperatures,
\[
  \Delta(\Lambda_{N},\eta^{N,0},\beta) \geq \exp(-\lambda(\delta)
    N^{\frac{1}{2} + \delta}) \quad \text{for all }
    \delta > 0.
\]
Presuming $N^{-2}$ is the correct rate for the ``plus'' boundary condition,
Theorem \ref{T:main}(i) shows that removing as few as $O(\log N)$ plus spins from
the corners of the box dramatically shrinks the spectral gap.

For $k \approx cN$ with $0 < c < 1$, Theorem \ref{T:main} shows that the spectral
gap decreases at least exponentially fast in $N$.  Schonmann \cite{Sc94} showed
that for all $\eta$, all $N$ and all $\beta > 0$, for some $C = C(\beta)$,
\[
  \Delta(\Lambda_{N},\eta,\beta) \geq \frac{C}{N}e^{-4\beta N},
\]
so the gap can never decrease faster than exponentially in $N$.

Our proof of Theorem \ref{T:main} will use the method suggested by
(\ref{E:indicators}):  we find an event $D$ with $\partial_{in} D$ much
smaller than $D$.
This event $D$ is a variant of the event that none of the four strips of ``+''
spins in $\partial \Lambda_{N}$ is connected to any of the other
strips by a path of ``+'' spins.

\section{Preliminaries} \label{S:prelims}
Throughout the paper, $c_{0}, c_{1},...$ and $\epsilon_{1},
\epsilon_{2},...$ will be used to represent constants which depend only on the
temperature (or other parameters) of the model.  We use $\epsilon_{i}$ for
constants which should be viewed as ``small.''

Our proof will make use of the Fortuin-Kasteleyn random cluster model, or briefly,
the \emph{FK model} (\cite{Fo1}, \cite{Fo2}, \cite{FK}; see also \cite{ACCN},
\cite{Gr96}) which is a graphical representation of the Potts model.  To discuss
this, we need some notation for bond configurations.  By a bond we mean an
unordered pair $\langle xy \rangle$ of adjacent sites in
$\mathbb{Z}^{2}$.  When convenient 
we view bonds as being open line segments in the plane; this should be
clear from the context.  Define the sets of bonds
\[
  \mB(\Lambda) = \{\langle xy \rangle: x,y \in \Lambda \}, \quad
  \omB(\Lambda) = \{\langle xy \rangle: x \in \Lambda \text{ or } y \in \Lambda \}.
\]
For general $A \subset \mathbb{R}^{2}$, we write $\mB(A)$ for $\mB(A \cap \ZZ)$.
Let $\Omega_{\Lambda} = \{0,1\}^{\omB(\Lambda)}$. A
\emph{bond configuration} is an element $\omega \in \Omega_{\Lambda}$; when
convenient we alternatively view $\omega$ as a subset of $\omB(\Lambda)$ or as a
subgraph of
$\graph$.  Bonds $e$ with $\omega_{e} = 1$ are \emph{open} in $\omega$;
those with $\omega_{e} = 0$ are \emph{closed}. Let $C(\omega)$ denote the
number of open clusters in $\omega$ which do not intersect $\partial \Lambda$. 
For $p \in [0,1]$ and $q > 0$, the FK model 
$P_{\Lambda,w}^{p,q}$ on $\graph$
with parameters $(p,q)$ and wired boundary condition is defined by the weights
\begin{equation} \label{E:FKweight}
  W(\omega) = p^{|\omega|}(1 - p)^{|\omB(\Lambda)| - |\omega|}q^{C(\omega)}
\end{equation}
Here $|\omega|$ means the number of open bonds in $\omega$.  More generally, given
$\rho \in \{0,1\}^{\omB(\Lambda)^{c}}$ we define $(\omega \rho)$ to be the bond
configuration on the full lattice which coincides with $\omega$ on $\omB(\Lambda)$
and with $\rho$ on $\omB(\Lambda)^{c}$.  Let $C(\omega \mid \rho)$ be the number
of open clusters of $(\omega \rho)$ which intersect $\Lambda$.  The FK model 
$P_{\Lambda,\rho}^{p,q}$ with
bond boundary condition $\rho$ is given by the weights in (\ref{E:FKweight}) with
$C(\omega)$ replaced by $C(\omega \mid \rho)$.  Alternately, given $\eta \in
\{-1,0,1\}^{\partial \Lambda}$ define
\begin{align} \label{E:Devent}
  V(\Lambda,\eta) = \bigl\{\omega \in \{0,1\}^{\omB(\Lambda)}: \eta_{x}
    &= \eta_{y} \text{ for every } x,y \in \partial \Lambda \text{ for which }
    x \lra y \text{ in } \omega, \\
  &\omega_{e} = 0 \text{ for all } e \in \{\langle xy \rangle:
    x \in \Lambda, y \in \partial\Lambda, \eta_{y} = 0\} \bigr\}. \notag
\end{align}
Here $x \lra y$ means there is a path of open bonds connecting $x$
to $y$.  The FK model $P^{p,q}_{\Lambda,\eta}$ with
site boundary condition $\eta$ is given by the weights in (\ref{E:FKweight}),
multiplied by $\delta_{V(\Lambda,\eta)}(\omega)$.  Taking $\eta_{x} = 0$ for all
$x$ gives the FK model with free boundary condition; we denote it
$P^{p,q}_{\Lambda,f}$.  For a summary of basic
properties of the FK model, see \cite{Gr96}.  In particular, since we are in two
dimensions, for $p \neq \sqrt{q}/(1 + \sqrt{q})$ there is a unique
translation-invariant infinite-volume FK measure on $\mB(\ZZ)$, which can be
obtained as the limit of $P^{p,q}_{\Lambda,w}$ as $\Lambda \nearrow \ZZ$; we denote
this measure $P^{p,q}$. For $q \geq 1$, the FK model has the FKG property.  For 
\[
  \alpha(p,q) = \frac{p}{p + q(1-p)}
\]
and $ e \in \mB(\ZZ)$, we have
\[
  P^{p,q}(\omega_{e} = 1 \mid \omega_{b}, b \neq e ) \geq \alpha(p,q) \quad
    \text{ for every } (\omega_{b}, b \neq e).
\]
Changing a single bond in the boundary condition changes the value of
$C(\omega \mid \rho)$ by at most 1.  It follows easily that for boundary conditions
$\rho, \rho^{\prime}$ differing at only one bond, we have
\begin{equation} \label{E:bondeffect}
  P^{p,q}_{\Lambda,\rho} \leq qP^{p,q}_{\Lambda,\rho^{\prime}}.
\end{equation}

As shown in \cite{ES}, for $\beta$ given by $p = 1 -
e^{-\beta}$, a configuration of the Ising model on $\Lambda$ with boundary
condition $\eta$ at inverse temperature $\beta$ can be obtained from a
configuration
$\omega$ of the FK model at $(p,2)$ with site boundary condition $\eta$, by
choosing a label for each cluster of
$\omega$ independently and uniformly from $\{-1,1\}$; this
\emph{cluster-labeling construction} yields a joint site-bond configuration for
which the sites are an Ising model and the bonds are an FK model.  When the
parameters are related in this way, we call the Ising and FK models
\emph{corresponding}.  Alternately, if one selects an Ising configuration $\sigma$
and does independent percolation at density $p$ on the set of bonds
\[
  \{\langle xy \rangle \in \omB(\Lambda): \sigma_{x} = \sigma_{y} \},
\]
the resulting bond configuration is a realization of the corresponding FK model. 
We call this the \emph{percolation construction} of the FK model.
From Onsager's exact solution of the
Ising model (see \cite{MW}) and basic properties of the FK model (see \cite{Gr96})
the critical point $\beta_{c}$ of the Ising model and the percolation critical
point $p_{c}$ for the FK model with $q = 2$ are given by
\begin{equation} \label{E:critpts}
  1 - e^{-\beta_{c}} = p_{c} = \frac{\sqrt{2}}{1 + \sqrt{2}}.
\end{equation}

The \emph{dual lattice} $(\mathbb{Z}^{2})^{*}$ is $\mathbb{Z}^{2}$ shifted by
$(\tfrac{1}{2},\tfrac{1}{2})$; sites and bonds of this lattice are called
\emph{dual sites} and \emph{dual bonds}.  $x^{*}$ denotes $x +
(\tfrac{1}{2},\tfrac{1}{2})$.  When necessary for clarity,
bonds of
$\mathbb{Z}^{2}$ are called \emph{regular bonds}.  To each regular bond $e$ there
is associated a unique dual bond $e^{*}$ which is its perpendicular bisector.  
For $\mathcal{D} \subset \mB(\ZZ)$ we write $\mathcal{D}^{*}$ for $\{ e^{*}: e \in
\mathcal{D} \}$.  For $\Lambda \subset (\ZZ)^{*},\ \partial \Lambda$ is defined
as for $\Lambda \subset \ZZ$, but using adjacency in the dual lattice.  The
dual bond $e^{*}$ is defined to be open precisely when $e$ is closed, so that for
each bond configuration $\omega$ on
$\mathbb{Z}^{2}$, there is unique dual configuration $\omega^{*}$ on
$(\mathbb{Z}^{2})^{*}$.  For $p \in [0,1]$ the value $p^{*}$
\emph{dual to} $p$ at level $q$ is given by
\[
  \frac{p}{q(1-p)} = \frac{1 - p^{*}}{p^{*}}.
\]
If the regular bonds are distributed as the infinite-volume FK model at
$(p,q)$ on $\mathbb{Z}^{2}$ with wired boundary condition, then the dual bonds
form the infinite-volume FK model at $(p^{*},q)$ on $(\mathbb{Z}^{2})^{*}$ with
free boundary condition (see \cite{Gr96}.) 

An Ising configuration $\sigma \in \Sigma_{\Lambda}$ determines a set of contours,
each consisting of dual bonds $e^{*} \in \omB(\Lambda)^{*}$ for which the
corresponding regular bond $e = \langle xy \rangle$ has $\sigma_{x} \neq
\sigma_{y}$.  In the joint Ising/FK configuration, therefore, 
contours consist entirely of open dual bonds.

Given sets $\Phi \subset \Lambda$ and a site configuration $\sigma \in
\Sigma_{\Lambda}$, we write $\sigma_{\Phi}$ for $\{ \sigma_{x}: x \in \Phi\ \}$
and let $\mathcal{F}_{\Phi}$ denote the $\sigma$-algebra generated by 
$\sigma_{\Phi}$.  Similarly for $\mathcal{D} \subset \omB(\Lambda)$ and a bond
configuration $\omega \in \Omega_{\Lambda}$, we write $\omega_{\mathcal{D}}$
for $\{\omega_{e}: e \in \mathcal{D}\}$ and let $\mathcal{G}_{\mathcal{D}}$ denote
the $\sigma$-algebra generated by $\omega_{\mathcal{D}}$.

An infinite-volume FK model $P^{p,q}$ (or other bond percolation model) is said to
have the
\emph{weak mixing property} if there exist $C,\lambda$ such that, given
finite sets
$\Phi \subset \Lambda$ and any two bond boundary conditions $\rho_{1}$ and
$\rho_{2}$ on $\omB(\Lambda)^{c}$, we have
\[
  \Var\bigl( P^{p,q}_{\Lambda,\rho_{1}}(\omega_{\mB(\Phi)} \in
    \cdot), P^{p,q}_{\Lambda,\rho_{2}}(\omega_{\mB(\Phi)} \in \cdot)\bigr)
    \leq C \sum_{x \in \Phi,y \notin \Lambda} e^{-\lambda |y - x|},
\]
Loosely this says that the maximum influence, on a
fixed region, of the boundary condition decays exponentially to 0 as the
boundary recedes to infinity.  Equivalently, for all events $A \in
\mathcal{G}_{\Lambda^{c}}$ and $B \in \mathcal{G}_{\Phi}$, 
\begin{equation} \label{E:weakmix}
  |P^{p,q}(B \mid A) - P^{p,q}(B)| \leq 
    C \sum_{x \in \Phi,y \notin \Lambda} e^{-\lambda |y - x|}.
\end{equation}
In contrast, $P^{p,q}$ is said to have the
\emph{ratio weak mixing property} if there exist $C,\lambda$ such that, given
finite sets
$\Phi \subset \Lambda$ and any two bond boundary conditions $\rho_{1}$ and
$\rho_{2}$ on $\omB(\Lambda)^{c}$, we have for all events $A \in
\mathcal{G}_{\Lambda^{c}}$ and $B \in \mathcal{G}_{\Phi}$,
\begin{equation} \label{E:rweakmix}
  \left| \frac{P^{p,q}(A \cap B)}{P^{p,q}(A)P^{p,q}(B)} - 1 \right|
    \leq C \sum_{x \in \Phi,y \notin \Lambda} e^{-\lambda |y - x|},
\end{equation}
whenever the right side of this inequality is at most 1.
Note that (\ref{E:rweakmix}) is much stronger than (\ref{E:weakmix}) for $A,B$ for
which the probabilities on the left side of (\ref{E:weakmix}) are much smaller
than the right side of (\ref{E:weakmix}).  Weak mixing for the Ising
model has a variety of useful consequences, particularly in two dimensions;
see \cite{MOS}.  It was shown in \cite{Al98} that for the FK model in two
dimensions, exponential decay of either the
connectivity (in infinite volume, with wired boundary) or the dual
connectivity (in infinite volume, with free boundary) implies ratio weak mixing. 
In particular, for
$q = 2$ and $p > p_{c}(2)$, exponential decay of dual connectivity follows from
the known properties (see \cite{MW}) of
Gibbs uniqueness and exponential decay of correlations for the Ising model
at inverse temperature $\beta^{*} < \beta_{c}(2)$ corresponding to $p^{*} <
p_{c}(2)$.  Thus we have the following.

\begin{lemma} \label{L:FKratiowm}
  Suppose $p > p_{c}(2)$.  Then the FK model $P^{p,2}$ on $\mB(\ZZ)$ has the
  ratio weak mixing property.
\end{lemma}

The following is an immediate consequence of the definition of
ratio weak mixing.

\begin{lemma} \label{L:ratiowm}
(\cite{Al97pwr})
Suppose that the FK model $P^{p,q}$ has the ratio weak mixing property.
There exists a constant $c_{1}$ as follows.
Suppose $r > 3$ and $\mathcal{D}, \mathcal{E} \subset \mB(\mathbb{Z}^{2})$
with $\diam(\mathcal{E}) \leq r$ and $d(\mathcal{D},\mathcal{E}) \geq
c_{1}\log r$.  Then for all $A \in \mathcal{G}_{\mathcal{D}}$ and
$B \in \mathcal{G}_{\mathcal{E}}$, we have
\[
  \frac{1}{2}P^{p,2}(A)P^{p,2}(B) \leq P^{p,2}(A \cap B) \leq 2P^{p,2}(A)P^{p,2}(B).
\]
\end{lemma}

We write $y \lrad z$ for the event that $y$ is connected to $z$ by a path of open
dual bonds.  For $q \geq 1$, $P^{p,q}$ has the FKG property (see \cite{Gr96}), so
$-\log P^{p,q}(0^{*} \lrad x^{*})$ is a subadditive function
of $x$, and therefore the limit
\begin{equation} \label{E:surftens}
  \tau(x) = 
  \lim_{n \to \infty} -\frac{1}{n}\log
  P^{p,q}(0^{*} \lrad (nx)^{*}),
\end{equation}
exists for $x \in \mathbb{Q}^{2}$, provided we take the limit
through values of $n$ for which $nx \in \ZZ$.  This definition 
extends to $\mathbb{R}^{2}$ by continuity (see \cite{Al97}); the 
resulting $\tau$ is a norm on $\mathbb{R}^{2}$, when the dual
connectivity decays exponentially (i.e. $\tau(x)$ is positive for
all $x \neq 0$, or equivalently by lattice symmetry, $\tau(x)$ is
positive for some $x \neq 0$; we abbreviate this by saying
$\tau$ \emph{is positive}.)
By standard subadditivity results,
\begin{equation} \label{E:conupr}
  P^{p,q}(0^{*} \lrad x^{*}) \leq e^{-\tau(x)}
  \quad \text{for all} \ x.
\end{equation}
In the opposite direction, it is known \cite{Al97pwr}
that if $\tau$ is positive (so
ratio weak mixing holds), then for some $\epsilon_{1}$ and $c_{2}$,
\begin{equation} \label{E:conlwr}
  P^{p,q}(0^{*} \lrad x^{*}) \geq
  \epsilon_{1}|x|^{-c_{2}}e^{-\tau(x)}
  \quad \text{for all} \ x \ne 0.
\end{equation}
It follows from the fact that the surface tension $\tau$
is a norm on $\mathbb{R}^{2}$ with axis symmetry that,
letting $e_{i}$ denote the $i$th unit coordinate vector, we have
\begin{equation} \label{E:equivnorm}
  \frac{1}{\sqrt{2}}\tau(e_{1}) \leq
  \frac{\tau(x)}{|x|} \leq \sqrt{2}\tau(e_{1}) \quad
  \text{ for all } \ x \ne 0.
\end{equation}
  
A weakness of Lemma \ref{L:ratiowm} is that the locations
$\mathcal{D}, \mathcal{E}$ of the two events must be deterministic.  
The next lemma from \cite{Al00} applies only to a limited class of events
but allows the locations to be partially random.
For $\mathcal{C} \subset \mathcal{D}
\subset \mathcal{B}(\ZZ)$ we say an event $A 
\subset \{0,1\}^{\mathcal{D}}$ 
\emph{occurs on} $\mathcal{C}$ (or on $\mathcal{C}^{*}$)
in $\omega \in \{0,1\}^{\mathcal{D}}$ 
if $\omega^{\prime} \in A$  for every
$\omega^{\prime}\in \{0,1\}^{\mathcal{D}}$  satisfying 
$\omega^{\prime}_{e} = \omega_{e}$
for all $e \in \mathcal{C}$.  
For a possibly random set $\mathcal{F}(\omega)$ we
say $A$ \emph{occurs only on} $\mathcal{F}$ (or equivalently,
on $\mathcal{F}^{*}$) if $\omega \in A$ implies
$A$ occurs on $\mathcal{F}(\omega)$ in $\omega$.
We say events $A$ and $B$ 
\emph{occur at separation} $r$ in $\omega$ if there exist 
$\mathcal{C}, \mathcal{E} \subset \mathcal{D}$ with 
$d(\mathcal{C},\mathcal{E}) \geq r$ such
that $A$ occurs on $\mathcal{C}$ and $B$ occurs on
$\mathcal{E}$ in $\omega$.  
Let $A \circ_{r} B$ denote the event that $A$ and $B$ 
occur at separation $r$.
Let $\mathcal{D}^{r} = \{e \in \mathcal{B}(\ZZ): 
d(e,\mathcal{D}) \leq r\}$.

For $x$ a (regular or dual) site, we write $C_{x} = C_{x}(\omega)$ for the
(regular or dual) cluster of
$x$ in the bond configuration $\omega$.

\begin{lemma} \label{L:decouple}
  (\cite{Al00})
  Let $P^{p,q}$ be an FK model on $\mathcal{B}(\ZZ)$, with $\tau$ positive
  and $q \geq 1$, satisfying
  the ratio weak mixing property.  
  There exist constants $c_{i}, \epsilon_{i}$ as follows.
  Let $\mathcal{D} \subset \mathcal{B}(\ZZ), x \in 
  (\ZZ)^{*}$ and 
  $r > c_{3} \log |\mathcal{D}|$, and let
  $A, B$ be
  events such that $A$ occurs only on $C_{x}$ and
  $B \in \mathcal{G}_{\mathcal{D}}$.  Then 
  \begin{equation} \label{E:decouple}
    P^{p,q}(A \circ_{r} B) \leq
    (1 + c_{4}e^{-\epsilon_{2}r}) P^{p,q}(A)P^{p,q}(B).
  \end{equation}
\end{lemma}

Let $\overline{xy}$ denote the line through $x$ and $y$.  Let 
\[
  H_{a}^{+} = \{ (x_{1},x_{2}) \in \mathbb{R}^{2}: x_{2} \geq a \}
\]
and let $H_{xy}$ denote the closed halfspace bounded by $\overline{xy}$ which is
to the right as one moves from $x$ to $y$.  

\begin{lemma} \label{L:halfspace}
  (\cite{Al97pwr})
  Let $P^{p,q}$ be an FK model on $\mathcal{B}(\ZZ)$, with $\tau$ positive
  and $q \geq 1$, satisfying
  the ratio weak mixing property.  There exist
  $\epsilon_{3}, c_{5}$ such that for all $x \neq y \in \mathbb{R}^{2}$
  and all dual sites $u, v \in H_{xy}$, 
  \[
    P(u \lrad v \text{ via a path in } 
    H_{xy}) \geq \frac{\epsilon_{3}}{|v-u|^{c_{5}}}e^{-\tau(v-u)}.
  \]
\end{lemma}

For the $\tau$ of (\ref{E:surftens}), the next lemma is an easy consequence
of the sharp triangle inequality satisfied by $\tau$ (see \cite{Io}), which is
obtained from the exact solution of the Ising model on $\ZZ$.  We will not use
this method here, however, to help make it apparent that our results are not
specific to the Ising model.

For $u \in \mathbb{R}^{2}$ let
$D_{u}^{+}$ and $D_{u}^{-}$ denote the diagonal lines through $u$ and $u + (1,1)$
and through $u$ and $u + (1,-1)$, respectively.  Let $T_{N,k}^{1}$ 
denote the triangle with vertices $(-k,-N), (k,-N)$ and
$(0,-(N-k))$.  Note the base
of $T_{N,k}^{1}$ is in the bottom side of $\Lambda_{N}$ and the other two
sides are parallel to the diagonals.  Let
$T_{N,k}^{i}, i = 2,3,4$, be the corresponding triangles (obtained by
rotation) with bases in the left, top and right sides of $\Lambda_{N}$,
respectively.

\begin{lemma} \label{L:normprop}
  Suppose $\tau$ is a norm on the plane which has axis and diagonal symmetry.  
  There exists a constant $\epsilon_{4}$ as follows.  Let 
  $0 < k < k+2m < N$ and $x = (x_{1},-N - \frac{1}{2}), y = (y_{1},-N
  - \frac{1}{2})$ with $x_{1} \leq -k,
  y_{1} \geq k$, and $z = (z_{1},z_{2}) \in H_{-N-\frac{1}{2}}^{+} \backslash
  T_{N+\frac{1}{2},k+2m+\frac{1}{2}}^{1}$.  Then
  \begin{equation} \label{E:normprop}
    \tau(z-x) + \tau(y-z) \geq 2k\tau(e_{1}) + \epsilon_{4}m.
  \end{equation}
\end{lemma}
\begin{proof}
If $x_{1} \leq -k-m$ or $y_{1} \geq k + m$ then $\tau(y - x) \geq
\tau((2k+m)e_{1})$ and (\ref{E:normprop}) follows easily.  Hence assume $|x_{1}|,
y_{1} \in [k,k+m)$.  We may also assume $z_{1} \geq 0$.  

If $z_{1} > y_{1}$ then from symmetry and convexity we have $\tau(z - x) \geq
\tau(y-x)$ and $\tau(y-z) \geq \tau(me_{1}/2)$, and (\ref{E:normprop}) follows
easily.  Hence we assume $0 \leq z_{1} \leq y_{1}$.  

If $z$ is above $D_{x}^{+}$, let $z^{\prime}$ be the reflection of $z$ across
$D_{x}^{+}$.  Then $z^{\prime} \notin T_{N,k+2m}^{1}$ and $\tau(z^{\prime}
- x) = \tau(z - x)$, and
it follows easily from symmetry and convexity that $\tau(y - z) \geq \tau(y -
z^{\prime})$.  Thus it is sufficient to prove (\ref{E:normprop}) for
$z^{\prime}$.  Hence we may assume $z$ is on or below $D_{x}^{+}$.

Now let $u = (u_{1},u_{2})$ be the reflection of $z$ across $D_{y}^{-}$, and let
$v = (v_{1},v_{2})$ be the point where $\overline{xz}$ intersects $D_{y}^{-}$.  By
the above assumptions on $z$ and simple geometry, we have $x_{1} \leq u_{1} \leq
v_{1}$ and $x_{2} \leq u_{2} \leq v_{2}$.  Using symmetry and convexity we
therefore obtain
\[
  \tau(z-x) = \tau(z-v) + \tau(v-x) \geq \tau(z-v) + \tau(u-x).
\]
Since $\tau(y-z) = \tau(y-u)$, it follows that
\begin{align}
  \tau(z-x) + \tau(y-z) &\geq \tau(z-v) + \tau(u-x) + \tau(y-u) \notag \\
  &\geq \epsilon_{4}m + \tau(y-x), \notag 
\end{align}
as desired.
\end{proof}

For $x, y \in (\ZZ)^{*}, r > 0$ and $G \subset \mathbb{R}^{2}$,
we say there is an \emph{r-near dual connection} from $x$ to $y$ 
in $G$ if for some $u, v \in (\ZZ)^{*}$ with $d(u,v) \leq r$,
there are open dual paths from $x$ to $u$ and from $y$ to $v$
in $G$.  Let $N(x,y,r,G)$ denote the event that such an $r$-near dual
connection exists.
The following result is from \cite{Al97pwr}.

\begin{lemma}\label{L:rnear}
  Let $P^{p,q}$ be an FK model on $\mathcal{B}(\ZZ)$, with $q \geq 1$, for
  which $\tau$ is positive.  There exist $c_{i}$ such that if 
  $|x| > 1$ and $r \geq c_{6}\log |x|$ then
  \[
    P^{p,q}(N(0,x,r,\mathbb{R}^{2})) \leq e^{-\tau(x) + c_{7}r}.
  \]
\end{lemma}

The next lemma shows that a dual connection via a site for which the triangle
inequality is strict by an amount $t > 0$ has an excess cost proportional to $t$.

\begin{lemma} \label{L:triangle}
  Let $P^{p,q}$ be an FK model on $\mathcal{B}(\ZZ)$, with $q \geq 1$, for which
  $\tau$ is positive.  There
  exists $c_{8}$ as follows.  Suppose $x, y, z \in
  (\ZZ)^{*}$ and $t \geq c_{8} \log |y-x|$ satisfy $|y-x| > 1$ and
  \[
    \tau(z-x) + \tau(y-z) \geq \tau(y-x) + t.
  \]
  Then
  \[
    P^{p,q}(x \lrad z \lrad y) \leq e^{-\tau(y-x) - \frac{1}{20}t}.
  \]
\end{lemma}
\begin{proof}
By Lemma \ref{L:FKratiowm}, $P^{p,q}$ has the ratio weak mixing property.  Let
$B = B_{\tau}(z,3(\tau(y-x) + t))$.  Then provided $c_{8}$ is large,
\begin{equation} \label{E:triangle1}
  P^{p,q}(z \lrad B^{c}) \leq c_{9}(\tau(y-x) + t)e^{-\tau(y-x) - t}
    \leq \frac{1}{3}e^{-\tau(y-x) - \frac{1}{2}t}.
\end{equation}
Thus we need only consider paths inside $B$.  Let
$\epsilon_{5} > 0$ to be specified
and consider $\omega \in [x \lrad
z \lrad y] \cap N(x,y,\epsilon_{5}t,B \backslash B_{\tau}(z,t/5))^{c}$. 
For such $\omega$, there exist $z^{\prime}, z^{\prime\prime} \in \partial(
B_{\tau}(z,t/5) \cap
(\ZZ)^{*})$ and paths $x \lrad z^{\prime}, y \lrad z^{\prime\prime}$ in $B$ 
occuring at separation $\epsilon_{5}t$.  Now
\[
  \tau(z^{\prime} - x) + \tau(y - z^{\prime\prime}) \geq \tau(x-x) + \tau(y-z)
    - 2\left(\frac{t}{5} + \tau(e_{1})\right) \geq \tau(y-x) + \frac{1}{2}t,
\]
so using Lemma \ref{L:decouple}, provided $c_{8}$ is large,
\begin{align} \label{E:triangle2}
  P^{p,q}&\bigl( [x \lrad z \lrad y \text{ in } B] \cap N(x,y,\epsilon_{5}t,B
    \backslash
    B_{\tau}(z,t/5))^{c} \bigr) \\   
  &\leq \sum_{z^{\prime},z^{\prime\prime}} 2 P^{p,q}(x \lrad z^{\prime}) \
    P^{p,q}(z^{\prime\prime} \lra y) \notag \\
  &\leq 2 |\partial(
    B_{\tau}(z,t/5) \cap (\ZZ)^{*})|^{2} e^{-\tau(y-x) - t/2}
    \notag \\
  &\leq \frac{1}{3}e^{-\tau(y-x) - t/4}. \notag
\end{align}
Next, provided $\epsilon_{5}$ is small, an application of Lemma \ref{L:ratiowm}
gives
\begin{align} \label{E:triangle3}
  P^{p,q}&\bigl( [x \lrad z \lrad y \text{ in } B] \cap N(x,y,\epsilon_{5}t,B
    \backslash
    B_{\tau}(z,t/5)) \bigr) \\
  &\leq P^{p,q}\bigl( [z \lrad B_{\tau}(z,t/10)^{c}] \cap 
    N(x,y,\epsilon_{5}t,B \backslash B_{\tau}(z,t/5)) \bigr) \notag \\
  &\leq 2 P^{p,q}\bigl( z \lrad B_{\tau}(z,t/10)^{c} \bigr)
    P^{p,q}\bigl( N(x,y,\epsilon_{5}t,B \backslash B_{\tau}(z,t/5))
    \bigr) \notag \\
  &\leq 2c_{10}te^{-t/10}e^{-\tau(y-x) + c_{7}\epsilon_{5}t} \notag \\
  &\leq \frac{1}{3}e^{-\tau(y-x) - t/20}. \notag
\end{align}
Together, (\ref{E:triangle1}), (\ref{E:triangle2}) and (\ref{E:triangle3})
complete the proof.
\end{proof}

\section{Proof of Theorem \ref{T:main}$(i)$} \label{S:mainproof0}
Let
\[
  \oL = \Lambda \cup \partial \Lambda, \quad \Lambda \subset \ZZ.
\]
By a \emph{plus path} in a site configuration $\sigma$ we mean a lattice path on
which all sites $z$ have $\sigma_{z} = 1$; \emph{minus paths} are defined
analogously.  We write $x \lrap y$ ($x \lram y$) for the event that $x$ is
connected to $y$ by a plus (minus) path.  For $\Phi \subset \oL_{N}$, the
\emph{cluster} of
$\Phi$ in a bond configuration $\omega \in \{0,1\}^{\omB(\Lambda_{N})}$ is the
set
\[
  C(\Phi,\omega) = \{ x \in \Lambda: \text{ in } \omega,
    x \lra \Phi \text{ in } \omB(\Lambda_{N}) \}.
\]
The \emph{plus cluster} of $\Phi$ in a
site configuration $\sigma \in \Sigma_{\Lambda_{N}}$ is the set
\[
  C_{+}(\Phi,\sigma) = 
    \{ x \in \Lambda: \text{ in } \sigma,
    x \lrap \Phi \text{ in } \omB(\Lambda_{N}) \}.
\]
If $\sigma_{x} = -1$, then of course $C_{+}(x,\sigma)$ is empty.
The \emph{minus cluster} $C_{-}(\Phi,\sigma)$ is defined analogously.

For $\Phi \subset \ZZ$ we define
\[
  Q(\Phi) = \cup_{x \in \Phi} \ \left(x +
    \left[-\frac{1}{2},\frac{1}{2}\right]^{2}\right).
\]
Here $x + [-\tfrac{1}{2},\tfrac{1}{2}]^{2}$ denotes the translation of 
$[-\tfrac{1}{2},\tfrac{1}{2}]^{2}$ by $x$.  For $x \in \mathbb{R}^{2}$ and $r > 0$
we let $B(x,r)$ and $B_{\tau}(x,r)$ be the closed Euclidean ball and $\tau$-ball,
respectively, of radius $r$ about $x$.

Fix $\beta > \beta_{c}, N \geq 1$ and $K \log N \leq k \leq N$, with $K$ to be
specified later.  Let $p = 1 - e^{-\beta}$ and $j = (N - k)/4$.  
We assume $k$ and $j$ are integers; the modifications otherwise are trivial.
Let $\mathbb{P}^{\epsilon}_{N,k}$ denote the
joint site/bond distribution obtained using the percolation construction of the FK
model, for which the site marginal distribution is
$\mu_{\Lambda_{N},\eta^{k,\epsilon}}^{\beta}$ and the bond marginal distribution
is $P_{\Lambda_{N},\eta^{k,\epsilon}}^{p,2}$.  (To avoid ambiguous notation we
write $\mathbb{P}^{-}_{N,k}$ for $\mathbb{P}^{-1}_{N,k}$.)  We write
$(\sigma,\omega)$ for a generic joint configuration in 
$\Sigma_{\Lambda_{N}} \times
\{0,1\}^{\omB(\Lambda_{N})}$.  We call $(\sigma,\omega)$ \emph{allowable} (under
$\eta^{k,\epsilon}$) if
$\mathbb{P}^{\epsilon}_{N,k}((\sigma,\omega)) > 0$.
Define the strip of sites $\Gamma_{1}$ by
$\Gamma_{1} = ([-k,k] \times \{-N-1\}) \cap \ZZ$, and let $\Gamma_{i}, 
i = 2,3,4,$ be the
corresponding strips of sites, obtained by rotation, in the left, top and right
sides of $\partial \tilde{\Lambda}_{N+1}$, respectively.  (We will refer to the
side corresponding to subscript $i$ as the $i$th side of $\tilde{\Lambda}_{n}$,
for general $n$.)  For
$i = 1,2,3$ let
$\Gamma_{i,i+1}$ be the set of sites in $\partial \Lambda_{N}$ which are between
$\Gamma_{i}$ and $\Gamma_{i+1}$, in the obvious sense, and let 
$\Gamma_{4,5} = \Gamma_{4,1} = \Gamma_{0,1}$ be
the set of sites in $\partial \Lambda_{N}$ which are between $\Gamma_{4}$ and
$\Gamma_{1}$.  We also include the appropriate ``corner site''
as an element of $\Gamma_{i,i+1}$, e.g. $(-N-1,-N-1) \in \Gamma_{1,2}$.  Let 
\[
  \Gamma^{*}_{i,j} = \{ x \in (\ZZ)^{*} \cap \partial\tilde{\Lambda}_{N +
    \tfrac{1}{2}}: 
    x \text{ is a corner of } Q(y) \text{ for some } y \in \Gamma_{i,j} \}.
\]
Let $D_{i}$ be the event that there is no plus-path in $\sigma$ in
$\omB(\Lambda_{N})$ from
$\Gamma_{i}$ to $(T_{N,k+3j}^{i})^{c}$, and $D = \cap_{i = 1}^{4} D_{i}$. 
In the FK model, an event closely related to $D_{i}$ is
\[
  E_{i} = \{ \omega \in \{0,1\}^{\omB(\Lambda_{N})}:  \Gamma^{*}_{i-1,i} 
    \lrad \Gamma^{*}_{i,i+1} \text{ in }
    T^{i}_{N+\frac{1}{2},k+j+\frac{1}{2}} \}.
\]

We begin with 
a lower bound on the probability of $E_{i}$.  The main point is that 
restricting the path to lie in $T^{i}_{N+\frac{1}{2},k+j+\frac{1}{2}}$ does not
excessively alter the probability of an open dual path from $\Gamma^{*}_{i-1,i}$ to
$\Gamma^{*}_{i,i+1}$.

\begin{lemma} \label{L:lowerbound}
  Let $p > p_{c}(2)$.  
  There exist $c_{i}, \epsilon_{i}$ such that for
  $N, k \geq 1$ with $c_{11} \log N \leq N - k \leq N$, and $E_{i}$ as
  above,
  \[
    P_{\Lambda_{N},\eta^{k,0}}^{p,2}(E_{i}) \geq \frac{\epsilon_{6}}{k^{c_{13}}}
      e^{-2k\tau(e_{1})}.
  \]
\end{lemma}
\begin{proof}
We may assume $i = 1$ and $k \geq c_{12}$, with $c_{12}$ to be specified. 
Let $m = c_{14} \log k$ and $n =
\lfloor c_{15} \log k \rfloor$, where $c_{14} > c_{15}$ are to be specified
and $\lfloor \cdot \rfloor$ denotes the integer part.
Provided $c_{11}$ is large (depending on $c_{14}$), we have $m \leq j$.  Let
$x^{\prime} = (-k-\tfrac{1}{2},-N-\tfrac{1}{2}), y^{\prime} =
(k+\tfrac{1}{2},-N-\tfrac{1}{2}), x = x^{\prime} + (n,n), y = y^{\prime} + (-n,n)$
and let 
\[
  E_{1}^{\prime} = \{ \omega \in \{0,1\}^{\omB(\Lambda)}:  x
    \lrad y \text{ in }
    T^{1}_{N,k+2m} \cap H_{-N+n}^{+} \}.
\]
By Lemma \ref{L:halfspace}, for some $\epsilon_{7}$,
\begin{equation} \label{E:halfsp}
  P^{p,2}(x \lrad y \text{ in }
    H_{-N+n}^{+}) \geq \frac{\epsilon_{7}}{k^{c_{5}}}e^{-2k\tau(e_{1})}.
\end{equation}
Note this probability is for the infinite-volume limit.
For each dual site $z \in H_{-N+n}^{+} \backslash T^{1}_{N,k+2m}$, let
$F_{z}$ denote the event that there exist open dual paths $x^{\prime} \lrad z \lrad
y^{\prime}$.  By (\ref{E:halfsp}), 
\begin{equation} \label{E:infvolbd}
  P^{p,2}(E_{1}^{\prime}) \geq \frac{\epsilon_{7}}{k^{c_{5}}}e^{-2k\tau(e_{1})}
    - P^{p,2}\left( \cup_{z \in (\ZZ)^{*} \cap H_{-N+n}^{+} \backslash 
    T^{1}_{N,k+2m}} F_{z}
    \right).
\end{equation}
We wish to show that the second term on the right side of (\ref{E:infvolbd}) is at
most half of the first term on the right side.  Let $\Theta = (\ZZ)^{*} 
\cap B_{\tau}(x,3k\tau(e_{1})) \cap H_{-N+n}^{+}
\backslash T^{1}_{N,k+2m}$.  We have
\begin{equation} \label{E:decomp}
  P^{p,2}\left( \cup_{z \in (\ZZ)^{*} \cap H_{-N+n}^{+} \backslash 
    T^{1}_{N,k+2m}} F_{z} \right) 
    \leq P^{p,2}\bigl( x \lrad B_{\tau}(x,3k\tau(e_{1}))^{c} \bigr) +
    \sum_{z \in \Theta} P^{p,2}(F_{z}). 
\end{equation}
By (\ref{E:conupr}), provided $c_{12}$ is large,
\begin{equation} \label{E:longpath}
  P^{p,2}\bigl( x \lrad B_{\tau}(x,3k\tau(e_{1}))^{c} \bigr) 
    \leq c_{16}ke^{-3k\tau(e_{1})} \leq
    \frac{1}{4}\frac{\epsilon_{7}}{k^{c_{5}}}e^{-2k\tau(e_{1})}.
\end{equation}
Fix $z \in \Theta$.  We decompose the event $F_{z}$ according to
whether there is a $(c_{6} \log 5k)$-near dual connection from $x$ to $y$ in
$B(z,3c_{17} \log k)^{c}$, where $c_{6}$ is from Lemma \ref{L:rnear} and
$c_{17}$ is to be specified.  
Let $\Psi_{z} = B(z,c_{17} \log k) \cap (\ZZ)^{*}$.  Using (\ref{E:conupr})
and Lemma \ref{L:ratiowm}, provided $c_{17}$ and $c_{12}$ are large we obtain
\begin{align} \label{E:yesrnear}
  P^{p,2}&\bigl( F_{z} \cap N(x,y,c_{6} \log 5k,B(z,2c_{17} \log k)^{c}) \bigr) \\
  &\leq P^{p,2}\bigl( [z \lrad \partial\Psi_{z}]
    \cap N(x,y,c_{6} \log 5k,B(z,2c_{17} \log k)^{c}) \bigr) \notag \\
  &\leq 2 P^{p,2}\bigl( z \lrad \partial\Psi_{z} \bigr)
    P^{p,2}\bigl( N(x,y,c_{6} \log 5k,B(z,2c_{17} \log k)^{c}) \bigr) \notag \\
  &\leq 2 c_{18}ke^{-\frac{1}{2}c_{17}\tau(e_{1}) \log k} e^{-\tau(y - x) +
    c_{7}c_{6} \log 5k} \notag \\
  &\leq e^{-\frac{1}{4}c_{17}\tau(e_{1}) \log k - 2(k-n)\tau(e_{1})}, \notag
\end{align}
where $c_{7}$ is from Lemma \ref{L:rnear}.
Since $|\Theta| \leq c_{19}k^{2}$,
provided we choose $c_{17}$ large enough (depending on $c_{15}$) this gives
\begin{equation} \label{E:yesrnear2}
  \sum_{z \in \Theta}
    P^{p,2}\bigl( F_{z} \cap N(x,y,c_{6} \log 5k,B(z,2c_{17} \log k)^{c}) \bigr)
    \leq \frac{1}{8}\frac{\epsilon_{7}}{k^{c_{5}}}e^{-2k\tau(e_{1})}.
\end{equation}
Next, let $r = c_{6} \log 5k$ and for $z \in \Theta$ let $\Psi_{z}^{\prime} =
B(z,2c_{17} \log k) \cap (\ZZ)^{*}$.  We have
\begin{equation} \label{E:uniondecoup}
  F_{z} \cap N(x,y,c_{6} \log 5k,B(z,2c_{17} \log k)^{c})^{c} \subset
    \cup_{u,v \in \partial \Psi_{z}^{\prime}} \bigl( [x \lra u] \circ_{r} [y \lra
    v] \bigr).
\end{equation}
Now for $u,v \in \partial \Psi_{z}^{\prime}$,
\[
  \tau(u-x) \geq \tau(z-x) - \tau(u-z) \geq \tau(z-x) - 3c_{17}\tau(e_{1}) \log k
\]
and similarly
\[
  \tau(y-v) \geq \tau(y-v) - 3c_{17}\tau(e_{1}) \log k.
\]
Hence provided $c_{14}$ is large enough (depending on $c_{17}$), we obtain using
Lemma \ref{L:normprop} that
\[
  \tau(u-x) + \tau(y-v) \geq 2(k-n)\tau(e_{1}) + \epsilon_{4}m - 6c_{17}\tau(e_{1})\log k
    \geq 2(k-n)\tau(e_{1}) + \frac{\epsilon_{4}}{2}m.
\]
Combining this with (\ref{E:uniondecoup}), Lemma \ref{L:decouple} and
(\ref{E:conupr}), provided $c_{14}$ and $c_{14}/c_{15}$ are large we get
\begin{align} \label{E:nornear}
  P^{p,2}\bigl( F_{z} \cap N(x,y,c_{6} \log 5k,B(z,2c_{17} \log k)^{c})^{c} \bigr)
    &\leq \sum_{u,v \in \partial \Psi_{z}^{\prime}} 2P^{p,2}(x \lrad u) P^{p,2}(y
    \lrad v) \\
  &\leq |\partial \Psi_{z}^{\prime}|^{2} e^{-2(k-n)\tau(e_{1}) - \epsilon_{4}m/2}
    \notag \\
  &\leq e^{-2k\tau(e_{1}) - \epsilon_{4}m/4}, \notag
\end{align}
and then
\begin{equation} \label{E:nornear2}
  \sum_{z \in \Theta} P^{p,2}\bigl( F_{z} \cap N(x,y,c_{6} 
    \log 5k,B(z,2c_{17} \log k)^{c})^{c} \bigr) 
    \leq c_{19}k^{2} e^{-2k\tau(e_{1}) - \epsilon_{4}m/4} 
    \leq \frac{1}{8}\frac{\epsilon_{7}}{k^{c_{5}}}e^{-2k\tau(e_{1})}. 
\end{equation}
Combining (\ref{E:infvolbd}), (\ref{E:decomp}), (\ref{E:yesrnear2}) and
(\ref{E:nornear2}) we obtain
\[
  P^{p,2}(E_{1}^{\prime}) \geq 
    \frac{1}{2}\frac{\epsilon_{7}}{k^{c_{5}}}e^{-2k\tau(e_{1})}.
\]
Then from Lemma \ref{L:ratiowm}, provided $c_{15}$ is large,
\[
  P^{p,2}_{\Lambda_{N},\eta^{k,0}}(E_{1}^{\prime}) \geq 
    \frac{1}{4}\frac{\epsilon_{7}}{k^{c_{5}}}e^{-2k\tau(e_{1})}.
\]
Let $\gamma_{x}$ and $\gamma_{y}$ be dual paths of (minimal) length $2n$ from $x$
to $x^{\prime}$ and from $y$ to $y^{\prime}$, respectively, in
$T^{1}_{N,k+2m}$.  Let $E_{1}^{\prime\prime}$ denote the event that all
dual bonds in $\gamma_{x}$ and $\gamma_{y}$ are open.  From the FKG inequality,
\[
  P^{p,2}_{\Lambda_{N},\eta^{k,0}}(E_{1}) \geq
    P^{p,2}_{\Lambda_{N},\eta^{k,0}}(E_{1}^{\prime} \cap E_{1}^{\prime\prime})
    \geq P^{p,2}_{\Lambda_{N},\eta^{k,0}}(E_{1}^{\prime})
    P^{p,2}_{\Lambda_{N},\eta^{k,0}}(E_{1}^{\prime\prime})
    \geq \frac{1}{4}\frac{\epsilon_{7}}{k^{c_{5}}}e^{-2k\tau(e_{1})} 
    \alpha(p,2)^{4n}
\]
and the lemma follows.
\end{proof}

For $\omega \in E_{i},\ \partial Q(C(\Gamma_{i},\omega))$ includes a unique open
dual path $\gamma_{i}(\omega)$ in $T^{i}_{N+\frac{1}{2},k+j+\frac{1}{2}}$
from $\Gamma^{*}_{i-1,i}$ to $\Gamma^{*}_{i,i+1}$.
This path is ``closer to $\Gamma_{i}$'' than any other open dual path in
$\tilde{\Lambda}_{N + \frac{1}{2}}$ from $\Gamma^{*}_{i-1,i}$ to
$\Gamma^{*}_{i,i+1}$.  Further, for fixed $\nu$ the event $[\gamma_{i} = \nu]$
depends only on the bond/dual bond configuration in the closed region, which we
denote $I(\nu)$, between $\nu$ and the side of $\partial \tilde{\Lambda}_{N +
\frac{1}{2}}$ to which $\nu$ is attached.

Let $E = \cap_{i = 1}^{4} E_{i}$, and suppose $\omega \in E$.  Let
\[
  R(\omega) = \tilde{\Lambda}_{N +
    \frac{1}{2}} \backslash \cup_{i = 1}^{4} I(\gamma_{i}(\omega)).
\]
Note $0 \in R(\omega)$, and (under boundary condition $\eta^{k,0}$) all the dual
bonds forming $\partial R(\omega)$ are open in $\omega$.  The latter means that
for fixed $U$, conditionally on $R(\omega) = U$ the configuration on $\mB(U)$ is
the FK model with free boundary condition.  

Let 
\[
  Y_{N} = \tilde{\Lambda}_{N+\frac{1}{2}} 
    \backslash \cup_{i=1}^{4} T^{i}_{N+\frac{1}{2},k+j+\frac{1}{2}}, \quad
  Y_{N}^{\prime} = \{ x \in Y_{N}: d(x,\partial Y_{N}) \geq j\}.
\]
Let $0 < h < j$ to be specified, let $w_{12} = (k+2j,-N+h)$, let $\lambda_{12}$ be
the vertical line from $w_{12}$ down to $\partial \tilde{\Lambda}_{N}$ at
$(k+2j,-N)$ and let $\chi_{12}$ be the
vertical line from $w_{12}$ up to the diagonal $D_{0}^{-}$ at $(k+2j,-k-2j)$.
Using axis symmetry we obtain 7 more corresponding points $w_{ij}$ and
paths $\chi_{ij},\lambda_{ij}$, for $i = 1,2,3,4$ and $j = 1,2$, 
with $w_{ij}$ at distance $h$
from side $i$ of $\tilde{\Lambda}_{N}$.  

We want to show that with high probability, there are no open dual paths 
starting from
$Y_{N}^{\prime}$, or from near $\chi_{ij}$, which reach $\partial Y_{N}$.
Let $y_{ijl1},y_{ijl2}$ be the endpoints of the dual bond which is dual to the
$l$th bond of $\chi_{ij}$.  
Then
\[
  d(y_{ijlm},\partial Y_{N}) \geq h +
    \frac{l}{\sqrt{2}} \quad \text{for all } i,j,l,m.
\]
For $x \in Y_{N}$ let
\[
  G_{x} = B\left( x,\frac{1}{2}d(x,\partial Y_{N}) \right) \cap (\ZZ)^{*},
\]
and define
\[
  \Psi_{N} = \bigl( Y_{N}^{\prime} \cap (\ZZ)^{*} \bigr) \cup 
    \{ y_{ijlm}:  1 \leq i \leq 4; j = 1,2; 1 \leq l \leq N-h-k-2j; m = 1,2 \}.
\]
Suppose $U \supset Y_{N}$.  Then $d(G_{x},\partial U) \geq j/2$ for all $x \in
Y_{N}^{\prime}$.  Hence
using Lemma \ref{L:ratiowm}, provided $K$ and $h$ are large enough we get
\begin{align} \label{E:bigcluster}
  P^{p,2}_{\Lambda_{N},\eta^{k,0}}&(x \lrad \partial G_{x} \text{ for some }
    x \in \Psi_{N} \mid \omega \in E, R(\omega) = U) \\
  &=P^{p,2}_{U \cap \ZZ,f}(x \lrad \partial G_{x} \text{ for some }
    x \in \Psi_{N}) \notag \\
  &\leq \sum_{x \in Y_{N}^{\prime} \cap (\ZZ)^{*}}
    P^{p,2}_{U \cap \ZZ,f}(x \lrad \partial G_{x}) +  \sum_{i,j,l,m} 
    P^{p,2}_{U \cap \ZZ,f}(y_{ijlm} \lrad \partial G_{y_{ijlm}}) \notag \\
  &\leq \sum_{x \in Y_{N}^{\prime} \cap (\ZZ)^{*}}
    2 P^{p,2}(x \lrad \partial G_{x}) +  \sum_{i,j,l,m} 
    2 P^{p,2}(y_{ijlm} \lrad \partial G_{y_{ijlm}}) \notag \\
  &\leq |Y_{N}^{\prime} \cap (\ZZ)^{*}| e^{-j\tau(e_{1})/2}
    + c_{20} \sum_{l \geq 1} (h+l)e^{-\frac{1}{4}\tau(e_{1})(h+l)} \notag \\
  &\leq \frac{1}{2}. \notag
\end{align}
Let $F$ denote the event that $x \lrad \partial G_{x}$ for no $x \in \Psi_{N}$, and
all bonds in $\lambda_{ij}$ are open for all $i, j$.  If $\omega \in E \cap F$,
then there is an open circuit in $Y_{N}$ surrounding
$\Psi_{N}$ and for each $i$, a portion of this open circuit, together with
$\lambda_{i1}$ and $\lambda_{i2}$, forms an open path in $T_{N,k+3j}^{i}
\backslash T_{N,k+j}^{i}$ from a site adjacent to $\Gamma_{i-1,i}$ to a
site adjacent to $\Gamma_{i,i+1}$.  When this occurs (with $\omega \in E \cap F$),
we call this circuit together
with all $\lambda_{ij}$ a \emph{blocking pattern}.  Note the blocking pattern is
contained in $R(\omega)$.  We have using (\ref{E:bigcluster}) and the FKG
inequality that for all $U \supset Y_{N}$,
\begin{align} \label{E:blocking}
  P^{p,2}_{\Lambda_{N},\eta^{k,0}}&(\text{there is a blocking pattern in } 
    R(\omega)
    \mid \omega \in E, R(\omega) = U) \\
  &= P^{p,2}_{U \cap \ZZ,f}(\text{there is a blocking pattern in } U) \notag \\
  &\geq P^{p,2}_{U \cap \ZZ,f}(F) \notag \\
  &\geq \frac{1}{2}P^{p,2}_{U \cap \ZZ,f}(\text{all bonds in } \lambda_{ij} \text{
    are open for all } i,j) \notag \\
  &\geq \frac{1}{2}\alpha(p,2)^{8h}. \notag
\end{align}
But considering the cluster-labeling construction of the joint Ising/FK
configuration, we see that if the configuration $\omega \in E$ has a blocking
pattern  in $R(\omega)$ and
all sites $x$ in the blocking pattern have $\sigma_{x} = -1$ (which occurs with
probability 1/2, given such $\omega$), then
$\sigma \in D$. Thus from (\ref{E:blocking}), the FKG property and Lemma
\ref{L:lowerbound}, for some $\epsilon_{8}$,
\begin{equation} \label{E:Dlowerbound}
  \mu^{\beta}_{\Lambda_{N},\eta^{k,0}}(D) \geq \frac{1}{4}\alpha(p,2)^{8h}
    P^{p,2}_{\Lambda_{N},\eta^{k,0}}(E) \geq \frac{\epsilon_{8}}{k^{4c_{13}}}
    e^{-8k\tau(e_{1})}.
\end{equation}

We turn now to upper bounds on
$\mu^{\beta}_{\Lambda_{N},\eta^{k,0}}(\partial_{in}D)$.  
Analogously to $\gamma_{i}(\omega)$,
for $\sigma \in D$, $\partial Q(C_{+}(\Gamma_{i},\sigma))$ includes a unique open
dual path $\gamma_{i}^{+}(\sigma)$ in $T^{i}_{N+\frac{1}{2},k+3j+\frac{3}{2}}$ from
$\Gamma^{*}_{i,i+1}$ to $\Gamma^{*}_{i-1,i}$, for each $i$. For fixed $\nu$ the
event $[\gamma_{i}^{+} = \nu]$ depends only on the site configuration in
$\overline{I(\nu) \cap \ZZ}$.

Suppose $\sigma \in \partial_{in}D$ and $\sigma^{x} \notin D$.  Then for some $i$
we have $x \in
\partial C_{+}(\Gamma_{i},\sigma)$, and either $x \in
\partial T^{i}_{N,k+3j+1}$ or there is an open dual circuit 
$\gamma$ in $\omega$ outside
$I(\gamma_{i}^{+}(\sigma))$ which includes an edge of $Q(x)$ and surrounds 
some site outside $T^{i}_{N,k+3j}$.  We can choose $\gamma$ to be the outer
boundary of $Q(\Phi)$ for some plus-cluster $\Phi$ in $\sigma$.  In this case
we call $\gamma$ an \emph{appendable circuit attachable at} $x$.

According to Lemma
\ref{L:ratiowm}, we can choose a constant $c_{21}$ as follows.  Let $Z_{N}^{i,i+1}
= \{ x \in \mathbb{R}^{2}: d(x,\Gamma_{i,i+1}) \leq c_{21} \log N \}$ and
$Z_{N} = \cup_{i} Z_{N}^{i,i+1}$.  Let $V_{closed}$ be the event that all bonds in
$\{ \langle xy \rangle: y \in \Lambda_{N}, x \in \Gamma_{i,i+1} \text{ for some } i
\}$ are closed.  Let $V_{open}$ be the event that all bonds in $\mB(\ZZ)
\backslash \omB(\Lambda_{N})$ are open.
(Note the boundary condition $\eta^{k,0}$ conditions $\omega$ on
$V_{open} \cap V_{closed}$.)  
Then for all events $A \in \mathcal{G}_{\omB(\Lambda_{N}) \backslash
\mB(Z_{N})}$, we have
\begin{equation} \label{E:Veffect}
  \frac{1}{2}P^{p,2}(A) \leq P^{p,2}(A \mid V_{closed}) \leq 2 P^{p,2}(A).
\end{equation}
We therefore call $Z_{N}$ the
\emph{free-boundary influence region}.
In particular, provided $K$ is large (depending on $c_{21}$), using the FKG
inequality and (\ref{E:Veffect}) we have
\begin{align} \label{E:bdryinf}
  P^{p,2}_{\Lambda_{N},\eta^{k,0}}(Z_{N}^{i-1,i} \lrad Z_{N}^{i,i+1}) 
    &= P^{p,2}(Z_{N}^{i-1,i} \lrad Z_{N}^{i,i+1} \mid V_{open} \cap V_{closed}) \\
  &\leq P^{p,2}(Z_{N}^{i-1,i} \lrad Z_{N}^{i,i+1} \mid V_{closed}) \notag \\
  &\leq 2 P^{p,2}(Z_{N}^{i-1,i} \lrad Z_{N}^{i,i+1}) \notag \\
  &\leq c_{22}(\log N)e^{-2(k - c_{21} \log N)\tau(e_{1})} \notag \\
  &\leq N^{c_{23}} e^{-2k\tau(e_{1})}. \notag
\end{align}
Our main task is roughly to show, using Lemma \ref{L:normprop}, that the
probability for a connection $\Gamma^{*}_{i-1,i} \lrad \Gamma^{*}_{i,i+1}$ 
(specifically, part of $\gamma_{i}^{+}(\sigma)$) which
does not stay inside $T^{i}_{N+\frac{1}{2},k+3j+\frac{3}{2}}$ is smaller 
than the right side of (\ref{E:Dlowerbound}) by at least a factor of
$e^{-\epsilon_{9}j}$, for some $\epsilon_{9}$.
We must decompose the
event $\partial_{in}D$ into several pieces according to the geometry of
the sets $C_{+}(\Gamma_{i},\sigma)$ and $C(\Gamma_{i},\omega)$.  The most
difficult case is that of leakage along the (free) boundary, in which
$\gamma_{i}^{+}(\sigma)$ goes outside $T^{i}_{N+\frac{1}{2},k+3j+\frac{3}{2}}$ 
by way of $Z_{N}$.

We define one more special dual path as follows.  For $\Phi \subset
\Lambda_{N}$ let
\[
  \hat{C}(\Phi,\omega) = \{ x \in \Lambda_{N} \backslash Z_{N}: \text{in }
    \omega, x \lra \Phi \text{ in } \omB(\Lambda_{N}) \backslash
    \mB(Z_{N}) \}.
\]
If $(\sigma,\omega)$ is allowable and 
$\sigma \in D$, then $\partial Q(\hat{C}(\Gamma_{i},\omega))$ includes a unique
open dual path from $Z_{N}^{i,i+1}$ to $Z_{N}^{i-1,i}$ in $\omB(\Lambda_{N})
\backslash \mB(Z_{N})$; we denote this path $\hat{\gamma}_{i}(\omega)$.  We have
\[
  \hat{\gamma}_{i}(\omega) \subset I(\gamma_{i}(\omega)) \subset
    I(\gamma_{i}^{+}(\omega)) \subset T^{i}_{N+\frac{1}{2},k+3j+\frac{3}{2}}.
\]
Let $u_{i}(\omega)$ and $v_{i}(\omega)$ be the starting 
and ending sites, respectively, of
$\hat{\gamma}_{i}(\omega)$ in $Z_{N}^{i,i+1}$ and $Z_{N}^{i-1,i}$, respecitvely.
Also define
\[
  W_{N}^{i,i+1} = \{ x \in \mathbb{R}^{2}: d(x,\Gamma_{i,i+1}) \leq 
    2\epsilon_{10} j \}, \qquad W_{N} = \cup_{i} W_{N}^{i,i+1},
\]
where $\epsilon_{10}$ is to be specified.  Provided $K$ is large enough (depending
on $\epsilon_{10}$ and $c_{21}$), we have $Z_{N}^{i,i+1} \subset W_{N}^{i,i+1}$.

\emph{Case 1.}  Consider $\sigma \in 
\partial_{in} D$, and $\omega$ with $(\sigma,\omega)$
allowable, for which for some $i$ there exist (in order) dual sites $x,z,y \in
\hat{\gamma}_{i}(\omega)$ with 
\begin{equation} \label{E:normsum}
  \tau(y-x) \geq (2k - 2c_{21} \log N)\tau(e_{1}), \qquad
  \tau(z-x) + \tau(y-z) \geq (2k + 4\epsilon_{10}j)\tau(e_{1}).
\end{equation}
We let $A_{1}$ denote the set of $(\sigma,\omega)$ for which this occurs, let
$J_{i}$ denote the set of all $(x,y,z) \in
(T^{i}_{N+\frac{1}{2},k+3j+\frac{3}{2}} \cap (\ZZ)^{*})^{3}$ 
for which (\ref{E:normsum}) holds,
and let $J_{i}^{\prime}$ denote the set of all
$(x,y,z) \in J_{i}$ which also satisfy
\[
  \tau(y-x) \geq (2k + 2\epsilon_{10}j)\tau(e_{1}).
\]
As in (\ref{E:bdryinf}), provided K is large enough (depending on
$\epsilon_{10}$ and $c_{21}$), using (\ref{E:Veffect}), (\ref{E:bdryinf})
and Lemmas \ref{L:ratiowm} and \ref{L:triangle}
we get
\begin{align} \label{E:case1}
  \mathbb{P}^{0}(A_{1}) &\leq \sum_{i=1}^{4} \sum_{(x,y,z) \in J_{i}}
    P^{p,2}_{\Lambda_{N},\eta^{k,0}} \bigl( x \lrad z \lrad y \text{
    in } T^{i}_{N+\frac{1}{2},k+3j+\frac{3}{2}} \backslash \mB(Z_{N}); \\
  &\qquad \qquad \qquad Z_{N}^{l-1,l} \lrad Z_{N}^{l,l+1} \text{
    in } T^{l}_{N+\frac{1}{2},k+3j+\frac{3}{2}} \backslash \mB(Z_{N}) 
    \text{ for all } l \neq i \bigr) \notag \\
  &\leq \sum_{i=1}^{4} \sum_{(x,y,z) \in J_{i}}
    2 P^{p,2}\bigl( x \lrad z \lrad y \text{
    in } T^{i}_{N+\frac{1}{2},k+3j+\frac{3}{2}} \backslash \mB(Z_{N}); \notag \\
  &\qquad \qquad \qquad Z_{N}^{l-1,l} \lrad Z_{N}^{l,l+1} \text{
    in } T^{l}_{N+\frac{1}{2},k+3j+\frac{3}{2}} \backslash \mB(Z_{N}) 
    \text{ for all } l \neq i \bigr) \notag \\
  &\leq \sum_{i=1}^{4} \sum_{(x,y,z) \in J_{i}}
    16 P^{p,2}\bigl( x \lrad z \lrad y \text{
    in } T^{i}_{N+\frac{1}{2},k+3j+\frac{3}{2}} \backslash \mB(Z_{N}) \bigr)
    \notag \\
  &\qquad \qquad \qquad \quad \cdot \prod_{l \neq i} 
    P^{p,2}\bigl(Z_{N}^{l-1,l} \lrad Z_{N}^{l,l+1} \text{
    in } T^{l}_{N+\frac{1}{2},k+3j+\frac{3}{2}} \backslash \mB(Z_{N}) 
    \bigr) \notag \\
  &\leq 16 \left( \sum_{i=1}^{4} \sum_{(x,y,z) \in J_{i}^{\prime}} e^{-\tau(y-x)}
    + \sum_{i=1}^{4} \sum_{(x,y,z) \in J_{i} \backslash J_{i}^{\prime}}
    e^{-\tau(y-x) - \frac{1}{10}\epsilon_{10}j\tau(e_{1})} \right) \bigl(
    N^{c_{23}}e^{-2k\tau(e_{1})} \bigr)^{3} \notag \\
  &\leq 16 \left( 4|J_{1}^{\prime}| e^{-(2k + 2\epsilon_{10}j)\tau(e_{1})}
    + 4|J_{1} \backslash J_{1}^{\prime}| e^{-(2k +
    \frac{1}{10}\epsilon_{10}j - 2c_{21}\log N)\tau(e_{1})} \right) \bigl(
    N^{c_{23}}e^{-2k\tau(e_{1})} \bigr)^{3}\notag \\
  &\leq e^{-(8k + \frac{1}{20}\epsilon_{10}j)\tau(e_{1})}. \notag
\end{align}

\emph{Case 2}.  Let 
\[
  R_{N}^{i} = T^{i}_{N+\frac{1}{2},k+3j+\frac{3}{2}} \cap 
    (T^{i}_{N+\frac{1}{2},k+j+\frac{1}{2}} \cup W_{N}).
\]
We may think of $R_{N}^{i}$ as a ``triangle with feet.''
Let $A_{2}$ denote the set of all $(\sigma,\omega) \in \partial_{in} D \backslash
A_{1}$ for which, for some $i$, there exists a dual site $z \in
T^{i}_{N+\frac{1}{2},k+3j+\frac{3}{2}} \backslash R_{N}^{i}$ which is either in
$\gamma_{i}^{+}(\sigma)$ or in some appendable circuit attachable at some $x \in
\partial C_{+}(\Gamma_{i},\sigma)$. 
(In particular, this means that $z \lrad \partial
(B(z,\epsilon_{10}j) \cap (\ZZ)^{*})$ and 
$B(z,2\epsilon_{10}j) \subset \tilde{\Lambda}_{N+1}$.)  
Suppose $(\sigma,\omega) \in A_{2}$.  We claim that
$B(z,2\epsilon_{10}j)
\cap \hat{\gamma}_{i}(\omega) = \phi$.  For all $u \in
\Gamma_{i,i+1}^{*}$ and $v \in \Gamma_{i-1,i}^{*}$, by Lemma
\ref{L:normprop} we have
\[
  \tau(z - u) + \tau(v - z) \geq 2k\tau(e_{1}) + \epsilon_{4}j.
\]
Therefore for all $u^{\prime} \in Z_{N}^{i,i+1}, v^{\prime} \in Z_{N}^{i-1,i}$ and
$z^{\prime} \in B(z,2\epsilon_{10}j)$, provided $K$ is large (depending on
$c_{21}$) and $\epsilon_{10}$ is small (depending on $\epsilon_{4}$),
\[
  \tau(z^{\prime} - u^{\prime}) + \tau(v^{\prime} - z^{\prime}) \geq 2k\tau(e_{1})
    + \frac{1}{2}\epsilon_{4}j - 2c_{21}(\log N)\tau(e_{1})
    - 4\epsilon_{10}j\tau(e_{1}) \geq (2k + 4\epsilon_{10}j)\tau(e_{1}).
\]
Taking $u^{\prime} = u(\omega), v^{\prime} = v(\omega)$ and comparing to
(\ref{E:normsum}) we see that since $(\sigma,\omega) \notin A_{1}$, we cannot have
$z^{\prime} \in \hat{\gamma}_{i}(\omega)$, proving our claim.  Therefore as in
(\ref{E:case1}),
\begin{align} \label{E:case2}
  \mathbb{P}^{0}_{N,k}(A_{2}) &\leq \sum_{i=1}^{4} \sum_{z \in
    (T^{i}_{N+\frac{1}{2},k+3j+\frac{3}{2}} \backslash R_{N}^{i}) \cap (\ZZ)^{*}}
    P^{p,2}_{\Lambda_{N},\eta^{k,0}} \bigl(
    z \lrad B(z,\epsilon_{10}j)^{c} \text{ in } B(z,\epsilon_{10}j + 1), \\
  &\qquad \qquad \qquad \qquad \qquad \qquad
    Z_{N}^{i-1,i} \lrad Z_{N}^{i,i+1} \text{ in }
    T^{i}_{N+\frac{1}{2},k+3j+\frac{3}{2}} \cap B(z,2\epsilon_{10}j)^{c} \backslash
    \mB(Z_{N}), \notag \\
  &\qquad \qquad \qquad \qquad \qquad \qquad
    Z_{N}^{l-1,l} \lrad Z_{N}^{l,l+1} \text{ in }
    T^{l}_{N+\frac{1}{2},k+3j+\frac{3}{2}} \backslash \mB(Z_{N}) 
    \text{ for all } l \neq i\bigr) \notag \\
  &\leq \sum_{i=1}^{4} \sum_{z \in
    (T^{i}_{N+\frac{1}{2},k+3j+\frac{3}{2}} \backslash R_{N}^{i}) \cap (\ZZ)^{*}}
    2 P^{p,2} \bigl(
    z \lrad B(z,\epsilon_{10}j)^{c} \text{ in } B(z,\epsilon_{10}j + 1), \notag \\
  &\qquad \qquad \qquad \qquad \qquad \qquad
    Z_{N}^{i-1,i} \lrad Z_{N}^{i,i+1} \text{ in }
    T^{i}_{N+\frac{1}{2},k+3j+\frac{3}{2}} \cap B(z,2\epsilon_{10}j)^{c} \backslash
    \mB(Z_{N}), \notag \\
  &\qquad \qquad \qquad \qquad \qquad \qquad
    Z_{N}^{l-1,l} \lrad Z_{N}^{l,l+1} \text{ in }
    T^{l}_{N+\frac{1}{2},k+3j+\frac{3}{2}} \backslash \mB(Z_{N}) 
    \text{ for all } l \neq i\bigr) \notag \\
  &\leq 32 \sum_{i=1}^{4} \sum_{z \in
    (T^{i}_{N+\frac{1}{2},k+3j+\frac{3}{2}} \backslash R_{N}^{i}) \cap (\ZZ)^{*}}
    P^{p,2}(z \lrad B(z,\epsilon_{10}j)^{c}) \prod_{l=1}^{4} P^{p,2} \bigl(
    Z_{N}^{l-1,l} \lrad Z_{N}^{l,l+1} \bigr) \notag \\
  &\leq c_{24} |T^{1}_{N+\frac{1}{2},k+3j+\frac{3}{2}} \cap (\ZZ)^{*}| j
    e^{-\frac{1}{2}\epsilon_{10}j\tau(e_{1})}
    \bigl( N^{c_{23}}e^{-2k\tau(e_{1})} \bigr)^{4} \notag \\
  &\leq e^{-(8k + \frac{1}{4}\epsilon_{10}j)\tau(e_{1})}. \notag
\end{align}

\emph{Case 3}.  Let $A_{3} = \partial_{in} D \backslash
(A_{1} \cup A_{2})$, and suppose $(\sigma,\omega) \in A_{3}$.  In this case, plus
spins are ``leaking along the boundary,'' in the following sense: for some $i$, 
we have $\gamma_{i}^{+}(\sigma) \subset R_{N}^{i}$, and either
$\gamma_{i}^{+}(\sigma)$ or some appendable circuit $\gamma$ 
contains a dual site $z \in W_{N}$ at one of the ``toes'' of $R_{N}^{i}$, that is,
in the right or left side of 
$T^{i}_{N+\frac{1}{2},k+3j+\frac{3}{2}}$.  
Note that such a $\gamma$ necessarily has 
$\gamma \cap T^{i}_{N+\frac{1}{2},k+3j+\frac{3}{2}} \subset R_{N}^{i}$, since
$(\sigma,\omega)
\notin A_{2}$.  We will assume $z$ is in an appendable
circuit, attachable at some site $x$; the case of $z \in \gamma_{i}^{+}(\sigma)$
is similar but slightly simpler.  From symmetry, we may also assume 
$i = 1$ and $z$ is in the right side of
$T^{1}_{N+\frac{1}{2},k+3j+\frac{3}{2}}$.  We let $A_{3}^{\prime}(x,z)$
denote the event that $A_{3}$ occurs with a specified choice of $x,z$, with $i =
1$ and with
$z$ in an appendable circuit, and let $A_{3}^{\prime} = \cup_{x,z}
A_{3}^{\prime}(x,z)$.  We define $A_{3}^{\prime\prime}(z)$ and
$A_{3}^{\prime\prime}$ analogously for the case of $z \in \gamma_{i}^{+}(\sigma)$.

The Ising model has the following ``bounded energy'' property:
\begin{equation} \label{E:bddener}
  \mu_{\Lambda_{N},\eta^{k,0}}^{\beta}(\sigma_{y} = 1 \mid \sigma_{w}, w \neq y)
    \geq \frac{1}{1 + e^{8\beta}} \quad \text{for all } (\sigma_{w}, w \neq y).
\end{equation}

Given a site $x$, let $\omega^{x}$ denote the
configuration given by
\[
  \omega^{x}_{e} = 
  \begin{cases}
    0, &\text{if } x \text{ is an endpoint of } e; \\
    \omega_{e}, &\text{otherwise}.
  \end{cases}
\]
Note that if $(\sigma,\omega)$ is allowable, then so is $(\sigma^{x},\omega^{x})$.
Let
\[
  B^{trunc}_{m} = B \cap \{ (y_{1},y_{2}) \in \mathbb{R}^{2}: y_{1} < m \}
    \quad \text{for } m > 0, B \subset \mathbb{R}^{2},
\]
\[
  \hat{\Gamma}_{1,2} = \bigl( [k,k+2j) \times \{-N\} \bigr) \cap \ZZ,
\]
\[
  \psi_{1,2} = \{ k+2j \} \times [-N,-N + 2\epsilon_{10}j]
\]
and
\[
  \hat{C}_{-}(\hat{\Gamma}_{1,2},\sigma^{x}) = 
    \quad \{ y \in (\Lambda_{N})^{trunc}_{k+2j}:  \text{ in } \sigma^{x}, y
    \lram \hat{\Gamma}_{1,2} \text{ in }
    (\tilde{\Lambda}_{N})^{trunc}_{k+2j} \}. 
\]
For $J \subset \ZZ$ define the boundary condition $\eta^{k,0,J}$ by
\begin{equation} \label{E:bdrycond2}
  \eta^{k,0,J}_{u} = 
  \begin{cases}
    1, &\text{if } u \in J; \\
    \eta^{k,0}_{u}, &\text{otherwise}.
  \end{cases}
\end{equation}
Fix $x,z$ and suppose $(\sigma,\omega) \in A_{3}^{\prime}(x,z)$.
There is then a plus path in $\sigma^{x}$ from 
$\Gamma_{1}$ to $\psi_{1,2}$ in $R_{N}^{1}$, and 
$\hat{C}_{-}(\hat{\Gamma}_{1,2},\sigma^{x})$ is contained in 
the region between this
plus path and $\partial (\Lambda_{N})^{trunc}_{k+2j}$.  Since 
\[
  |\partial
  \hat{C}_{-}(\hat{\Gamma}_{1,2},\sigma^{x}) \cap \psi_{1,2}| \leq 2\epsilon_{10}j,
\]
using the ``bounded energy'' property (\ref{E:bddener}) of the Ising model we have
for fixed $J$
\begin{align} \label{E:addplus}
  \mathbb{P}^{0}_{N,k}&\bigl( (\sigma,\omega) \in A_{3}^{\prime}(x,z), 
    \hat{C}_{-}(\hat{\Gamma}_{1,2},\sigma^{x}) = J \bigr) \\
  &\leq \mathbb{P}^{0}_{N,k}\bigl( 
    \hat{C}_{-}(\hat{\Gamma}_{1,2},\sigma^{x}) = J,
    \text{ and in } \omega^{x}, \psi_{1,2} - \tfrac{1}{2}e_{1} \lrad
    Z_{N}^{4,1} \text{ in } \omB(\Lambda_{N} \backslash \overline{J})^{*} \cap
    (R_{N}^{1})^{trunc}_{k+2j-\frac{1}{2}} \notag \\
  &\qquad \qquad \qquad \text{and } Z_{N}^{l-1,l} \lrad Z_{N}^{l,l+1} \text{ in }
    T^{l}_{N+\frac{1}{2},k+j+\frac{1}{2}} \backslash \mB(Z_{N})
    \text{ for } l = 2,3,4 \bigr) \notag \\
  &\leq c_{25}^{2\epsilon_{10}j} \mathbb{P}^{0}_{N,k}\bigl( 
    \hat{C}_{-}(\hat{\Gamma}_{1,2},\sigma^{x}) = J,
    \text{ and in } \omega^{x}, \psi_{1,2} - \tfrac{1}{2}e_{1} \lrad
    Z_{N}^{4,1} \text{ in } \omB(\Lambda_{N} \backslash \overline{J})^{*} \cap
    (R_{N}^{1})^{trunc}_{k+2j-\frac{1}{2}} \notag \\
  &\qquad \qquad \qquad \text{and } Z_{N}^{l-1,l} \lrad Z_{N}^{l,l+1} \text{ in }
    T^{l}_{N+\frac{1}{2},k+j+\frac{1}{2}} \backslash \mB(Z_{N})
    \text{ for } l = 2,3,4 \notag \\
  &\qquad \qquad \qquad \text{and } \sigma_{y} = 1 \text{ for all }
    y \in \psi_{1,2} \cap \partial J
    \bigr) \notag \\
  &= c_{25}^{2\epsilon_{10}j} \mathbb{P}^{0}_{N,k}\bigl( 
    \text{in } \omega^{x}, \psi_{1,2} - \tfrac{1}{2}e_{1} \lrad
    Z_{N}^{4,1} \text{ in } \omB(\Lambda_{N} \backslash \overline{J})^{*} \cap
    (R_{N}^{1})^{trunc}_{k+2j-\frac{1}{2}} \notag \\
  &\qquad \qquad \qquad \text{and } Z_{N}^{l-1,l} \lrad Z_{N}^{l,l+1} \text{ in }
    T^{l}_{N+\frac{1}{2},k+j+\frac{1}{2}} \backslash \mB(Z_{N})
    \text{ for } l = 2,3,4 \mid \hat{C}_{-}(\hat{\Gamma}_{1,2},\sigma^{x}) =
    J \notag \\
  &\qquad \qquad \qquad \text{and } \sigma_{y} = 1 \text{ for all } y \in
    \psi_{12} \cap \partial J \bigr) \notag \\
  &\qquad \qquad \cdot \mathbb{P}^{0}_{N,k} \bigl( 
    \hat{C}_{-}(\hat{\Gamma}_{1,2},\sigma^{x}) =
    J \text{ and } \sigma_{y} = 1 \text{ for all } y \in \psi_{1,2} \cap \partial J
    \bigr). \notag
\end{align}
Here $\psi_{1,2} - \tfrac{1}{2}e_{1}$ means the translate of $\psi_{1,2}$ by
$-\tfrac{1}{2}e_{1}$.  When $\hat{C}_{-}(\hat{\Gamma}_{1,2},\sigma^{x}) = J$,
we have by definition $\sigma^{x}_{y} = 1$ for all $y \in \partial J \cap
(\Lambda_{N})^{trunc}_{k+2j}$ and for all $y \in \hat{\Gamma}_{1,2} \backslash J$,
so from  the Markov property of the Ising model, 
the conditioning on the right side of
(\ref{E:addplus}) is equivalent to conditioning on $\sigma^{x}_{y} = 1$ for all
$y \in \overline{J} \cup \hat{\Gamma}_{1,2}$.  Therefore the first probability on
the right side of (\ref{E:addplus}) is
\begin{align} \label{E:addplus2}
  &P^{p,2}_{\Lambda_{N} \backslash (\overline{J} \cup
    \hat{\Gamma}_{1,2}),\eta^{k,0,J}}
    \bigl( \text{in } \omega^{x}, \psi_{1,2} - \tfrac{1}{2}e_{1} \lrad
    Z_{N}^{4,1} \text{ in } \omB(\Lambda_{N} \backslash (\overline{J} \cup
    \hat{\Gamma}_{1,2}))^{*} \cap
    (R_{N}^{1})^{trunc}_{k+2j-\frac{1}{2}} \\
  &\qquad \qquad \qquad \qquad \text{and } 
    Z_{N}^{l-1,l} \lrad Z_{N}^{l,l+1} \text{ in }
    T^{l}_{N+\frac{1}{2},k+j+\frac{1}{2}} \backslash \mB(Z_{N})
    \text{ for } l = 2,3,4 \bigr). \notag 
\end{align}
Let
\[
  \hat{Z}_{N}^{1,2} = \{ x \in \mathbb{R}^{2}: d(x,\Gamma_{1,2} \backslash
    \hat{\Gamma}_{1,2}) \leq c_{21} \log N \}, \quad \hat{Z}_{N} = 
    \hat{Z}_{N}^{1,2} \cup Z_{N}^{2,3} \cup Z_{N}^{3,4} \cup Z_{N}^{4,1},
\]
where $c_{21}$ is from the definition of $Z_{N}$.
One effect of changing the measure from $P^{p,2}_{\Lambda_{N},\eta^{k,0}}$ to 
$P^{p,2}_{\Lambda_{N} \backslash (\overline{J} \cup
\hat{\Gamma}_{1,2}),\eta^{k,0,J}}$ is to shrink the free-boundary influence region
from $Z_{N}$ to $\hat{Z}_{N}$.  More precisely, as in (\ref{E:case1}),
we have using Lemma \ref{L:ratiowm} that (\ref{E:addplus2}) is at most
\begin{align} \label{E:case3}
  P^{p,2}_{\Lambda_{N} \backslash (\overline{J} \cup
    \hat{\Gamma}_{1,2}),\eta^{k,0,J}} &\biggl( \text{in } \omega^{x}, 
    Z_{N}^{4,1} \lrad 
    \left( \psi_{1,2} - \frac{1}{2}e_{1} \right) \cup \hat{Z}_{N}^{1,2} \\
  &\qquad \qquad \text{ in } \omB(\Lambda_{N} \backslash (\overline{J} \cup
    \hat{\Gamma}_{1,2}))^{*} \cap
    (R_{N}^{1})^{trunc}_{k+2j-\frac{1}{2}} \backslash 
    \mB(\hat{Z}_{N}) \text{ and } \notag \\
  &\qquad \qquad Z_{N}^{l-1,l} \lrad Z_{N}^{l,l+1} \text{ in }
    T^{l}_{N+\frac{1}{2},k+j+\frac{1}{2}} \backslash \mB(\hat{Z}_{N})
    \text{ for } l = 2,3,4 \biggr) \notag \\
  &\leq 16 P^{p,2}\biggl(\text{in } \omega^{x}, Z_{N}^{4,1} \lrad 
    \left( \psi_{1,2} - \frac{1}{2}e_{1} \right) \cup \hat{Z}_{N}^{1,2} \notag \\
  &\qquad \qquad \qquad \text{ in } \omB(\Lambda_{N} \backslash (\overline{J} \cup
    \hat{\Gamma}_{1,2}))^{*} \cap
    (R_{N}^{1})^{trunc}_{k+2j-\frac{1}{2}} \backslash 
    \mB(\hat{Z}_{N}) \biggr) \notag \\
  &\qquad \cdot \prod_{l=2}^{4} P^{p,2}\bigl( Z_{N}^{l-1,l} 
    \lrad Z_{N}^{l,l+1} \text{ in }
    T^{l}_{N+\frac{1}{2},k+j+\frac{1}{2}} \backslash \mB(\hat{Z}_{N}) 
    \bigr). \notag 
\end{align}
To bound the first probability on the right side of (\ref{E:case3}), we observe
that, provided $K$ is large, if $u \in \hat{Z}_{N}^{4,1}$ and $v \in
(\psi_{1,2} - \frac{1}{2}e_{1}) \cup \hat{Z}_{N}^{1,2}$, then
\[
  \tau(v-u) \geq (2k + 2j - \tfrac{1}{2} - 2c_{21}
    \log N)\tau(e_{1}) \geq (2k + j)\tau(e_{1}).
\]
Further, we can replace $\omega^{x}$ with $\omega$ at the expense of at most a
constant factor.  Therefore as in (\ref{E:case1}), the right side of
(\ref{E:case3}) is at most
\[
  c_{26}je^{-(2k+j)\tau(e_{1})}\bigr( N^{c_{23}}e^{-2k\tau(e_{1})}
    \bigr)^{3} \leq e^{-(8k + \frac{1}{2}j)\tau(e_{1})}.
\]
Plugging this into (\ref{E:addplus}), provided $\epsilon_{10} < 1/8$ and $K$ is
large we obtain
\begin{align} \label{E:A3prime}
  \mathbb{P}^{0}_{N,k}(A_{3}^{\prime}) &\leq \sum_{x,z,J} c_{25}^{2\epsilon_{10}j}
    e^{-(8k+\frac{1}{2}j)\tau(e_{1})} \mathbb{P}^{0}_{N,k} \bigl( 
    \hat{C}_{-}(\hat{\Gamma}_{1,2},\sigma^{x}) = J \bigr) \\
  &\leq e^{-(8k+\frac{1}{4}j)\tau(e_{1})}. \notag
\end{align}
A similar proof gives the same bound for $\mathbb{P}^{0}_{N,k}
(A_{3}^{\prime\prime})$. 
Combining this with (\ref{E:case1}) and (\ref{E:case2}) gives
\begin{equation} \label{E:bdryDbound}
  \mu_{\Lambda_{N},\eta^{k,0}}^{\beta}(\partial_{in}D) 
    = \mathbb{P}^{0}_{N,k}(A_{1} \cup A_{2} \cup A_{3})
    \leq e^{-(8k+\epsilon_{11}j)\tau(e_{1})}.
\end{equation}

It follows easily from (\ref{E:bdryinf}) that
$\mu_{\Lambda_{N},\eta^{k,0}}^{\beta}(D) \leq 1/2$.  Combining this with
(\ref{E:Dlowerbound}), (\ref{E:bdryDbound}) and (\ref{E:indicators}) yields
\[
  \Delta(\Lambda_{N},\eta^{k,0},\beta) 
    \leq c_{27}|\Lambda_{N}|k^{4c_{13}}e^{-\epsilon_{11}j}
    \leq e^{-\epsilon_{11}j/2} = e^{-\epsilon_{11}(N-k)/8},
\]
which proves Theorem \ref{T:main} for $\epsilon = 0$.

\section{Proof of Theorem \ref{T:main}$(ii)$}

The FK measure corresponding to the boundary condition $\eta^{k,-1}$ is given (cf.
(\ref{E:Devent})) by 
\[
  P^{p,2}_{\Lambda_{N},\eta^{k,-1}} = P^{p,2}_{\Lambda_{N},w}(\cdot \mid 
    V(\Lambda_{N},\eta^{k,-1})),
\]
where $V(\Lambda_{N},\eta^{k,-1})$ is the event that there is no open dual path
from $\Gamma_{i}$ to $\Gamma_{j,j+1}$ for any $i \neq j$.  Our calculations,
however, are facilitated by using a different conditioning, as follows.  
Consider bond configurations on $\omB(\Lambda_{N+1})$.  Let
$U^{FK}$ denote the event that for all $i$, all bonds $\langle xy \rangle$ with $x,
y \in \Gamma_{i}$ are open, all bonds $\langle xy \rangle$ with $x, y \in
\Gamma_{i,i+1}$ are open, and all other bonds in 
$\omB(\Lambda_{N+1}) \backslash \omB(\Lambda_{N})$ are
closed.   The FK model $P^{p,2}_{\Lambda_{N+1},f}(\cdot \mid
U^{FK})$ corresponds to an Ising model $\mu_{\Lambda_{N+1},f}^{\beta}(\cdot \mid
U^{Ising})$, where $U^{Ising}$ is the event that for all $i$, all sites in
$\Gamma_{i}$ have the same spin, and all sites in $\Gamma_{i,i+1}$ have the same
spin.  Let $L$ denote the event that for all $i$, all sites in $\Gamma_{i}$
have spin 1 and all sites in $\Gamma_{i,i+1}$ have spin -1.  Then 
\[
  \mu_{\Lambda_{N+1},f}^{\beta}(\sigma_{\Lambda_{N}}\in \cdot \mid
    L) = \mu_{\Lambda_{N},\eta^{k,-1}}.
\]  
The measure $\mathbb{P}^{0}_{N+1,k}(\cdot \mid 
\omega \in U^{FK}) = \mathbb{P}^{0}_{N+1,k}(\cdot
\mid \sigma \in U^{Ising})$ gives the joint construction, coupling 
$P^{p,2}_{\Lambda_{N+1,k},f}(\cdot \mid U^{FK})$ and 
$\mu_{\Lambda_{N+1},f}^{\beta}(\cdot \mid U^{Ising})$.

For $\sigma \in L$, let
\[
  \mathcal{J} = \{ (i,j): 1 \leq i,j \leq 4, i < j \},
\]
\[
  \mathcal{A}_{N}(\sigma) = \{ (i,j)\in \mathcal{J}: \Gamma_{i} \lrap
    \Gamma_{j} \text{ in } \omB(\Lambda_{N}) \}.
\]
As motivation, note we expect that, roughly,
\begin{align} \label{E:heuristic}
  \mu_{\Lambda_{N+1},f}^{\beta}(\mathcal{A}_{N} = \mathcal{J} \mid
    L) \approx 1
    &\quad \text{if } 2k\tau(e_{1}) > (N-k)\tau(e_{1}+e_{2}), \\
  \mu_{\Lambda_{N+1},f}^{\beta}(\mathcal{A}_{N} = \phi \mid
    L) \approx 1
    &\quad \text{if } 2k\tau(e_{1}) < (N-k)\tau(e_{1}+e_{2}). \notag 
\end{align}
In the case $2k\tau(e_{1}) \geq (N-k)\tau(e_{1}+e_{2})$, 
we will bound the spectral gap using in (\ref{E:indicators}) the same event $D$ as
in Section \ref{S:mainproof0}, but in the opposite case we replace it with a
different event
$\hat{D} = \cap_{i = 1}^{4} \hat{D}_{i,i+1}$.  Here $\hat{D}_{1,2}$ is the event
that there is no minus-path in $\sigma$ in $\omB(\Lambda_{N})$ from $\Gamma_{1}$ to
$(S^{1,2}_{N,k/4})^{c}$, where $S^{1,2}_{N,m}$ is the square $[m,N+1] \times
[-N-1,-m]$, and $\hat{D}_{i,i+1}, S^{i,i+1}_{N,m}$ are the corresponding event and
square obtained by rotation, for $i = 2,3,4$.

Suppose first that $2k\tau(e_{1}) \geq (N-k)\tau(e_{1}+e_{2})$.
Let $x_{1,1} = (-k-\tfrac{1}{2},-N-\tfrac{1}{2})$ and $x_{1,2} =
(k+\tfrac{1}{2},-N-\tfrac{1}{2})$.  These dual sites are approximately at the ends
of $\Gamma_{1}$.  We define corresponding sites $x_{ij}$ for $i = 2,3,4$ and $j =
1,2$.  In place of the event $E_{i}$ of Lemma \ref{L:lowerbound}, we will use
\[
  \hat{E}_{i} = \{ \omega:  x_{1,1} \lrad x_{1,2} \text{ in } 
    T^{i}_{N+\frac{1}{2},k+j+\frac{1}{2}} \}.
\]
Lemma \ref{L:lowerbound} and its proof remain valid for $\hat{E}_{i}$ in place of
$E_{i}$, and the proof of the lower bound (\ref{E:Dlowerbound}) for
$\mu^{\beta}_{\Lambda_{N},\eta^{k,0}}(D)$ goes through with minimal changes to give
\begin{equation} \label{E:Dlowerbound2}
  \mu^{\beta}_{\Lambda_{N+1},f}(D \cap L \mid U^{Ising}) \geq
    \frac{\epsilon_{12}}{k^{c_{28}}}e^{-8k\tau(e_{1})}.
\end{equation}
The proof of (\ref{E:bdryDbound}) also goes through with minimal changes; in fact
Case 3 can be made simpler using the fact that the boundary regions
$\Gamma_{i,i+1}$ are each wired. (We will not do so here, since it is
unnecessary.)  The result is that
\begin{equation} \label{E:bdryDbound2}
  \mu^{\beta}_{\Lambda_{N+1},f}(\partial_{in} D \cap L \mid
    U^{Ising}) \leq e^{-(8k+\epsilon_{13}j)\tau(e_{1})}.
\end{equation}
Combining (\ref{E:Dlowerbound2}) and (\ref{E:bdryDbound2}) gives
\begin{align} \label{E:muratio}
  \frac{\mu^{\beta}_{\Lambda_{N},\eta^{k,-1}}(\partial_{in} D)}
    {\mu^{\beta}_{\Lambda_{N},\eta^{k,-1}}(D)}
    &= \frac{\mu^{\beta}_{\Lambda_{N+1},f}(\partial_{in} D \mid L)}
    {\mu^{\beta}_{\Lambda_{N+1},f}(D \mid L)} \\
  &= \frac{\mu^{\beta}_{\Lambda_{N+1},f}(\partial_{in} D \cap L
    \mid U^{Ising})}
    {\mu^{\beta}_{\Lambda_{N+1},f}(D \cap L \mid U^{Ising})}
    \notag \\
  &\leq \frac{\epsilon_{12}}{k^{c_{28}}}e^{-\epsilon_{13}j\tau(e_{1})}. \notag
\end{align}

In Section \ref{S:mainproof0} we easily obtained the lower bound
$\mu^{\beta}_{\Lambda_{N},\eta^{k,0}}(D^{c}) \geq 1/2$ to complete the proof. 
Here the situation is a little more complex.  A lower bound of the form
\begin{equation} \label{E:Dclowerbd2}
  \frac{\mu^{\beta}_{\Lambda_{N+1},f}(D^{c} \cap L \mid U^{Ising})}
    {\mu^{\beta}_{\Lambda_{N+1},f}(D \cap L \mid U^{Ising})}
    \geq \theta
\end{equation}
for some $\theta$ is equivalent to the statement
\begin{equation} \label{E:Dclowerbd}
  \mu^{\beta}_{\Lambda_{N},\eta^{k,-1}}(D^{c}) 
    = \mu^{\beta}_{\Lambda_{N+1},f}(D^{c} \mid L)
    \geq \frac{\theta}{1 + \theta}.
\end{equation}
Hence we consider bounds for the numerator and denominator of
(\ref{E:Dclowerbd2}).  Let $F_{i,i+1}$ denote the event that $x_{i,2} \lrad
x_{i+1,1}$ via a path in $S^{i,i+1}_{N,k/4}$, and $F = \cap_{i=1}^{4} F_{i,i+1}$. 
We have
\begin{align} \label{E:DcversusF}
  \mu^{\beta}_{\Lambda_{N+1},f}(D^{c} \cap L \mid U^{Ising})
    &\geq \mu^{\beta}_{\Lambda_{N+1},f}(\hat{D} 
    \cap L \mid U^{Ising}) \\
  &\geq P^{p,2}_{\Lambda_{N+1},f}(F \mid U^{FK}) \mathbb{P}^{0}_{N+1,k}(\hat{D}
    \cap L \mid F \cap U^{FK}). \notag 
\end{align}
It is straightforward to prove an analog of Lemma \ref{L:normprop} for 
$S^{1,2}_{N+\frac{1}{2},2m+\frac{1}{2}}$ in place of
$T^{1}_{N+\frac{1}{2},2m+\frac{1}{2}}$.  Therefore mimicking the proof of Lemma
\ref{L:lowerbound}, we obtain
\begin{equation} \label{E:Fibound}
  P^{p,2}_{\Lambda_{N+1},f}(F_{i,i+1} \mid U^{FK}) \geq \frac{\epsilon_{14}}
    {(N-k)^{c_{29}}}e^{-(N-k)\tau(e_{1}+e_{2})}.
\end{equation}
Then, analogously to (\ref{E:Dlowerbound}), from (\ref{E:DcversusF}),
\begin{equation} \label{E:Dclowerbd3}
  \mu^{\beta}_{\Lambda_{N+1},f}(D^{c} \cap L \mid U^{Ising})
    \geq \frac{\epsilon_{15}}{(N-k)^{c_{29}}}e^{-4(N-k)\tau(e_{1}+e_{2})}.
\end{equation}
Next we have, using (\ref{E:bondeffect}) and Lemma \ref{L:ratiowm},
\begin{align} \label{E:Dbound}
  \mu^{\beta}_{\Lambda_{N+1},f}(D \cap L \mid U^{Ising})
    &\leq P^{p,2}_{\Lambda_{N+1},f}(x_{i,2} \lrad x_{i+1,1} \text{ in }
    T^{i}_{N + \frac{1}{2},k+3j+\frac{3}{2}} \text{ for all } i \mid U^{FK}) \\
  &\leq 2^{8} P^{p,2}_{\Lambda_{N},w}(x_{i,2} \lrad x_{i+1,1} \text{ in }
    T^{i}_{N + \frac{1}{2},k+3j+\frac{3}{2}} \text{ for all } i) \notag \\
  &\leq 2^{8} P^{p,2}(x_{i,2} \lrad x_{i+1,1} \text{ in }
    T^{i}_{N + \frac{1}{2},k+3j+\frac{3}{2}} \text{ for all } i) \notag \\
  &\leq 2048e^{-8k\tau(e_{1})}. \notag
\end{align}
Since $2k\tau(e_{1}) \geq (N-k)\tau(e_{1}+e_{2})$, (\ref{E:Dclowerbd3}) and
(\ref{E:Dbound}) give
\[
  \frac{\mu^{\beta}_{\Lambda_{N+1},f}(D^{c} \cap L \mid U^{Ising})}
    {\mu^{\beta}_{\Lambda_{N+1},f}(D \cap L \mid U^{Ising})}
    \geq \frac{\epsilon_{15}}{8(N-k)^{c_{29}}}.
\]
With (\ref{E:Dclowerbd2}) and (\ref{E:Dclowerbd}), this shows
\[
  \mu^{\beta}_{\Lambda_{N},\eta^{k,-1}}(D^{c}) \geq 
    \frac{\epsilon_{15}}{16(N-k)^{c_{29}}},
\]
which with (\ref{E:muratio}) completes the proof of Theorem \ref{T:main} for
$\epsilon = -1$, as in Section \ref{S:mainproof0}.

The proof when $2k\tau(e_{1}) < (N-k)\tau(e_{1}+e_{2})$ is similar, with the
roles of $D$ and $\hat{D}$ interchanged, using squares $S^{i,i+1}_{\cdot,\cdot}$
in place of the triangles $T^{i}_{\cdot,\cdot}$.

\section{Acknowledgement}
The author would like to thank N. Yoshida for helpful conversations.


\begin{thebibliography}{99}

\bibitem{ACCN}  Aizenman, M., Chayes, J.T., Chayes, L., and Newman, 
C.M., \emph{Discontinuity of the magnetization in the} $1/|x - y|^{2}$
\emph{Ising and Potts models}, J. Stat. Phys. \textbf{50} (1988), 1-40.

\bibitem{Al97} Alexander, K.S., \emph{Approximation of
subadditive functions and rates of convergence in limiting
shape results}, Ann. Probab. \textbf{25} (1997), 30-55.

\bibitem{Al97pwr}  Alexander, K.S., \emph{Power-law
corrections to exponential decay of connectivities
and correlations in lattice models}, preprint (1997).

\bibitem{Al98} Alexander, K.S., \emph{On weak mixing
in lattice models}, Probab. Theory Rel. Fields \textbf{110}
(1998), 441-471.

\bibitem{Al00} Alexander, K.S., \emph{Cube-root boundary fluctuations for droplets
in random cluster models}, preprint (2000).

\bibitem{ES} Edwards, R.G. and Sokal, A.D., \emph{Generalization of the
Fortuin-Kasteleyn-Swendsen-Wang representation and Monte Carlo algorithm},
Phys. Rev. D \textbf{38} (1988), 2009-2012.

\bibitem{FH} Fisher, D.S. and Huse, D.A., \emph{Dynamics of 
droplet fluctuations in
pure and random Ising systems}, Phys. Rev. B \textbf{35} (1987), 6841-6846.

\bibitem{Fo1} Fortuin, C.M., \emph{On
the random cluster model. II. The percolation model}, Physica \textbf{58}
(1972), 393-418.

\bibitem{Fo2} Fortuin, C.M., \emph{On
the random cluster model. III. the simple random-cluster process}, Physica
\textbf{59} (1972), 545-570.

\bibitem{FK} Fortuin, C.M. and Kasteleyn, P.W., \emph{On
the random cluster model. I. Introduction and relation to other models},
Physica \textbf{57} (1972), 536-564.

\bibitem{Gr96} Grimmett, G.R., \emph{The stochastic
random-cluster process and uniqueness of random-cluster measures},
Ann. Probab. \textbf{23} (1995), 1461-1510.

\bibitem{GP} Grimmett, G.R. and Piza, M.S.T., \emph{Decay of correlations in
random-cluster models}, Commun. Math. Phys. \textbf{189} (1997), 465-480.

\bibitem{HY} Higuchi, Y. and Yoshida, N., \emph{Slow relaxation of} $2-D$
\emph{stochastic Ising models with random and nonrandom boundary conditions}, in:
\emph{New Trends in Stochastic Analysis} (K.D. Elworthy, S. Kusuoka and I.
Shigekawa, eds.), 153-167, World Scientific, Singapore (1997).

\bibitem{Io} Ioffe, D., \emph{Large deviations for the} $2D$ \emph{Ising model:
A lower bound without cluster expansions}, J. Stat. Phys. \textbf{74} (1994),
411-432.

\bibitem{Li85} Liggett, T.M., \emph{Interacting Particle Systems},
Springer-Verlag, New York (1985).

\bibitem{Ma} Martinelli, F.,\emph{On the two dimensional dynamical Ising model in
the phase coexistence region}, J. Stat. Phys. \textbf{76} (1994), 1179-1246.

\bibitem{MOS} Martinelli, F., Olivieri, E., and Schonmann, R.H., \emph{For 2-D
lattice spin systems weak mixing implies strong mixing}, Commun. Math. Phys.
\textbf{165} (1994), 33-47.

\bibitem{MW} McCoy, B.M. and Wu, T.T., \emph{The Two-Dimensional Ising Model},
Harvard University Press, Cambridge, USA (1973).

\bibitem{Sc94} Schonmann, R.H., \emph{Slow droplet-driven relaxation of the
stochastic Ising models in the vicinity of the phase coexistence region}, Commun.
Math. Phys. \textbf{161} (1994), 1-49.

\bibitem{Th} Thomas, L.E., \emph{Bound on the mass gap for finite volume
stochastic Ising models at low temperature}, Commun. Math. Phys. \textbf{126}
(1989), 1-11.

\end{thebibliography}
\end{document}